\documentclass[12pt]{amsart}

\usepackage{graphicx}
\usepackage[leqno]{amsmath}
\usepackage{amssymb}
\usepackage[all]{xy}
\usepackage{graphicx}
\usepackage[usenames,dvipsnames]{color}
\usepackage[normalem]{ulem}

\newtheorem{thm}{Theorem}[section]
\newtheorem{notation}[thm]{Notation}
\newtheorem{lemma}[thm]{Lemma}
\newtheorem{cor}[thm]{Corollary}
\newtheorem{prop}[thm]{Proposition}

\newtheorem{remark}[thm]{Remark}

\theoremstyle{definition}
\newtheorem{definition}[thm]{Definition}

\date{\today}

\begin{document}

\title[$C^{1,0}$ Foliation Theory]
{$C^{1,0}$ Foliation Theory}

\author[Kazez]{William H.  Kazez}
\address{Department of Mathematics, University of Georgia, Athens, GA 30602}
\email{will@math.uga.edu}

\author[Roberts]{Rachel Roberts}
\address{Department of Mathematics, Washington University, St.  Louis, MO 63130}
\email{roberts@math.wustl.edu}

\keywords{codimension one foliation, flow, measured foliation, holonomy, flow box decomposition, Denjoy blowup}

\thanks{This work was partially supported by grants from the Simons Foundation (\#244855 to William Kazez, \#317884 to Rachel Roberts)}

\subjclass[2010]{Primary 57M50}

\begin{abstract}
Transverse one dimensional foliations play an important role in the study of codimension one foliations.  In \cite{KR2}, the authors introduced the notion of flow box decomposition of a 3-manifold $M$.  This is a combinatorial decomposition of $M$ that reflects both the structure of a given codimension one foliation and that of a given transverse flow, and that is  amenable to inductive strategies.  

In this paper, flow box decompositions are used to extend some classical foliation results to foliations that are not $C^2$.  Enhancements of well-known results of Calegari on smoothing leaves, Dippolito on Denjoy blowup of leaves, and Tischler on approximations by fibrations are obtained.  The methods developed are not intrinsically 3-dimensional techniques, and should generalize to prove corresponding results for codimension one foliations in $n$-dimensional manifolds.
\end{abstract}

\maketitle

\section{Introduction}

Smoothness plays an important role in the theory of codimension one foliations of 3-manifolds.  Reeb constructed the first $C^\infty$ foliation on $S^3$ as the union of two foliated solid tori, or Reeb components \cite{Reeb}.  This work of Reeb, together with work of Alexander \cite{Alexander} and Wallace \cite{Wallace}, led to the proofs by Lickorish \cite{Lickorish} and Novikov and Zieschang \cite{Novikov} that any closed 3-manifold has a $C^\infty$ codimension one foliation.  On the other hand, Haefliger \cite{Haefliger} showed that no foliation of $S^3$ can be analytic.  This was greatly improved by Novikov \cite{Novikov} to show that any $C^2$ foliation of $S^3$ must have Reeb components, and these never exist in analytic foliations.

The qualitative nature of foliation theory and its impact on the ambient 3-manifold was considerably advanced by Thurston's introduction of the norm on the homology of a 3-manifold, and in particular the minimizing properties of leaves of taut, transversely oriented, $C^2$ foliations \cite{Thurston}.  

A foliation is {\sl taut} if closed smooth transversals to the foliation pass through every point of the manifold.  This is also known as {\sl everywhere taut} to distinguish it from the more familiar notion of {\sl smoothly taut} in which closed smooth transversals are only required to intersect every leaf of the foliation.  For a discussion of these and other notions of tautness, why they are different for $C^{k,0}$ foliations, the same for $C^{k,1}$ foliations, and interchangeable up to $C^0$ approximation and isotopy of foliations, see \cite{KR5}.  

Foliations as a tool for understanding problems in 3-dimensional topology came to the fore as a result of Gabai's constructions of both $C^{\infty}$ and   often less smooth, but  finite depth, taut foliations \cite{G1,G2,G3}.  The success of Gabai's applications of foliation theory led to many constructions of taut codimension one foliations.  Often these foliations are constructed using Denjoy blowup techniques that yield foliations that are only $C^{\infty,0}$; that is, leaves are smoothly immersed, but transversely, their tangent plane fields vary only continuously.

The impetus for our work starts with the Eliashberg-Thurston approximation theorem \cite{ET}.  They showed that a taut, co-oriented codimension one $C^2$ foliation of a 3-manifold can be $C^0$ approximated by a pair of symplectically fillable contact structures.  This allows non-trivial Heegaard-Floer invariants to be assigned to any manifold that supports a taut foliation \cite{OS}.  This is, consequently, one of the pillars of the conjectural relationship between $L$-spaces, taut foliations, and left orderability of the fundamental group.  For details, see, for example \cite{OS2} and \cite{BGW}.

In \cite{KR2} and \cite{KR3} we extended the Eliashberg-Thurston approximation theorem to the class of all $C^{1,0}$, co-oriented taut foliations, thereby extending its reach to manifolds carrying the new constructions of foliations mentioned above. Similar results can be found in \cite{bowden}.  In doing so, we found that many of the standard tools for working with foliations had either not been developed for foliations with lesser smoothness than originally intended, or had not been developed with an eye towards $C^0$ approximation theory in which it is often necessary to produce a new foliation while only moving the tangent planes of the original foliation slightly.

This paper includes enhancements to well-known results of Calegari \cite{calegari} on smoothing leaves, Dippolito \cite{Di} on Denjoy blowup of leaves, and Tischler \cite{Ti} on approximations by fibrations.  It is possible that some of our results can be obtained by ``reading between the lines'' of the original source.   However, as is well-known,  subtleties, sometimes fatal, arise when smooth objects are replaced by objects that are merely continuous.  (See, for example, \cite{KR5}.)  These are foundational results in foliation theory, and proofs of these theorems, in the generality in which they are used, do not exist in the literature.  An advantage to the flow box decomposition approach we use is that each of these results can be proved directly with a single inductive strategy.

The methods developed in this paper are not intrinsically 3-dimen\-sion\-al techniques, and we expect they can be adapted to prove corresponding results for codimension one foliations in $n$-dimensional manifolds.

Basic definitions (codimension one foliation, flow, $(\mathcal F,\Phi)$ compatible, $C^0$ close, and $C^0$ small) are given in \S~2.  In \S~3 we recall the definition of flow box decomposition, define regular neighborhood structure, and prove a sequence of useful local smoothing results.  The main result of \S~4 is a proof that any $C^{1,0}$ foliation is isotopic to a $C^0$ close $C^{\infty,0}$ foliation.  In \S~5 we prove that any $C^{1,0}$ measured foliation is isotopic to a $C^0$ close smooth measured foliation.  Basic facts from \cite{KR2,KR3} about holonomy neighborhoods are recalled in \S~6.  We give Dippolito's definition \cite{Di} of Denjoy blowup in \S~7 and prove that particularly nice, $C^0$ close, Denjoy blowups of a $C^{1,0}$ codimension one foliation always exist.  

Throughout this paper, unless stated otherwise, $M$ will denote a 3-manifold that is either smooth or smooth with corners.  When $\partial M\ne\emptyset$, it is often useful to think of $M$ as a sutured manifold, not necessarily orientable, in the sense of \cite{G1}.  Recall that any topological 3-manifold admits a smooth structure, unique up to diffeomorphism \cite{Moise0,Moise}.

\section{Codimension one foliations and transverse flows}

We begin by defining foliations in 3-manifolds with empty boundary.  Near the end of this section, we extend these definitions to 3-manifolds with nonempty boundary that are smooth or smooth with corners; namely, manifolds locally modelled by open sets in $[0,\infty)^3$.

\begin{definition}\label{folndefn1} Let $M$ be a smooth 3-manifold with empty boundary.  Let $k$ and $l$ be non-negative integers or infinity with $l \leq k$.  Both $C^k$ and $C^{k,l}$ {\sl codimension one foliations} $\mathcal F$ are decompositions of $M$ into a disjoint union of $C^k$ immersed connected surfaces, called the {\sl leaves} of $\mathcal F$, together with a collection of charts $U_i$ covering $M$, with $\phi_i:\mathbb R^2 \times \mathbb R \to U_i$ a homeomorphism, such that the preimage of each component of a leaf intersected with $U_i$ is a horizontal plane.  

The foliation $\mathcal F$ is $C^k$ if the charts $(U_i,\phi_i)$ can be chosen so that each $\phi_i$ is a $C^k$ diffeomorphism.  

The foliation $\mathcal F$ is $C^{k,l}$ if for all $i$ and $j$,

\begin{enumerate}

\item the derivatives $\partial_x^{\hskip .015 truein a}\partial_y^{\hskip .015 truein b}\partial_z^{\hskip .015 truein c}$, taken in any order, on the domain of each $\phi_i$ and each transition function $\phi_j^{-1}\phi_i$ are continuous for all $a + b \leq k$, and $c \leq l$, and

\item if $l \geq1$, $\phi_i$ is a $C^1$ diffeomorphism.

\end{enumerate}
\end{definition}

\begin{remark} The smoothness conditions on both the charts and the transition functions are to ensure that the smooth structure on the leaves is compatible with the smooth structure on $M$.  
\end{remark}

In particular, $T\mathcal F$ exists and is continuous if and only if $\mathcal F$ is $C^{1,0}$.  Also notice that $C^{k,l}$ foliations are $C^l$, but not conversely.

Two $C^{k,0}$ foliations $\mathcal F$ and $\mathcal G$ of $M$ are called {\sl $C^{k,0}$ equivalent} if there is a self-homeomorphism of $ M$ that maps the leaves of $\mathcal F$ to the leaves $\mathcal G$, and is $C^k$ when restricted to any leaf of $\mathcal F$.

We use the terms {\sl transverse, transversal,} and {\sl transversely} in the smooth sense; that is, they refer to smooth objects intersecting so that the associated tangent spaces intersect minimally.  In contrast, a curve is {\sl topologically transverse} to $\mathcal F$ if no nondegenerate subarc is isotopic, relative to its endpoints, into a leaf of $\mathcal F$.

Given a codimension one foliation $\mathcal F$, it is useful to fix a one dimensional foliation $\Phi$ transverse or topologically transverse to $\mathcal F$.  Such a $\Phi$ always exists and can be realized as the union of curves $\phi_p(t)$ of continuous local flows $\phi$.  When $\mathcal F$ is transversely oriented, $\Phi$ can be realized as the union of curves $\phi_p(t)$ of a global flow $\phi: M \times \mathbb R \to M$.  When $\mathcal F$ is $C^0$, this is proved in \cite{HH} (Theorems 1.1.2 and 1.3.2).  When $\mathcal F$ is $C^{1,0}$, $\Phi$ can be chosen to be smooth; in fact, in this case, it is elementary to see that $\Phi$ exists and consists of the integral curves of a smooth line field transverse to $T\mathcal F$.  See, for example, Lemma~5.1.1 of \cite{CC}.

\smallskip
\noindent {\bf Conventions:} Unless otherwise stated, throughout the rest of this paper, ``foliation'' will mean a codimension one transversely oriented foliation of a 3-manifold $M$.  Such $M$ will be assumed to be compact and oriented.  Since the foliations studied will be $C^{1,0}$ we can assume without loss of generality that a smooth transverse flow to the foliation is chosen.  To simplify the exposition, we will abuse language, and refer to a one dimensional foliation transverse to a codimension one foliation as a flow.
\smallskip

When a foliation $\mathcal F$ is understood, a submanifold of positive codimension in $M$ is called {\sl horizontal} if each component is a submanifold of a leaf of $\mathcal F$ and {\sl vertical} if each component is transverse to $\mathcal F$.  When both a foliation $\mathcal F$ and a transverse flow $\Phi$ in $M$ are understood, a submanifold of positive codimension in $M$ is called {\sl vertical} if and only if it can be expressed as a union of subsegments of the flow $\Phi$.  A codimension-0 submanifold $X$ of $M$ is called {\sl $(\mathcal F,\Phi)$-compatible} if its boundary is piecewise horizontal and vertical, and hence $\mathcal F$ and $\Phi$ restrict naturally to foliation and flow on $X$.  If $X$ is $(\mathcal F,\Phi)$-compatible, let $\partial_v X$ denote its vertical boundary, and let $\partial_h X$ denote its horizontal boundary.  

\begin{definition}
Suppose $X$ is an $(\mathcal F,\Phi)$-compatible submanifold of $M$, where possibly $X=M$.  An isotopy of $X$ which maps each flow segment of $\Phi|_X$ to itself is called a {\sl flow compatible}, or {\sl $\Phi$ compatible, isotopy}.  Note that a flow compatible isotopy of $X$ fixes $\partial_h X$ pointwise.
\end{definition}

By allowing the foliation atlas to include boundary charts, Definition~\ref{folndefn1} naturally extends to the case that $M$ has nonempty boundary that is either smooth or smooth with corners.  Smooth boundary components must either be a leaf of $\mathcal F$, and hence horizontal, or transverse to $\mathcal F$, and hence vertical.  A boundary component with corners must decompose along its corners into smooth subsurfaces, where if two subsurfaces share a corner, one is horizontal and one is vertical.  Such an $M$ is a sutured manifold, in the sense of \cite{G1}.  Thus, if $\partial M\ne\emptyset$ and we double $(M,\mathcal F)$ along $\partial_v M$, $D\mathcal F$ is a foliation of $DM$ with all components of $\partial (DM)$ leaves of $\mathcal F$.

We restrict attention to flows $\Phi$ that meet $\partial M$ in a similarly constrained way.  A flow is required to be either everywhere transverse or everywhere tangent to a smooth component of $\partial M$.  And if $(S,\gamma)$ is a boundary component with annular sutures $A(\gamma)$, a flow is required to be transverse to $R(\gamma)$ and tangent to $A(\gamma)$.  In particular, if $\Phi$ is a flow transverse to $\mathcal F$, it is possible to double $\Phi$ along $\partial_v M$ so that $M$ is a $(D\mathcal F,D\Phi)$-compatible submanifold of $DM$.

The terms $C^0$ close and $C^0$ small both refer to distances between tangent planes.  More formally, suppose that a metric, $d$, has been chosen on the set of continuous sections of the Grassmann bundle of 2-planes in $TM^3$.  Given a section, typically the tangent bundle of a foliation, $T\mathcal F$, we say that another section, $T\mathcal G$, can be found {\sl $C^0$ close} to $T\mathcal F$, if for all $\epsilon >0$ a $\mathcal G$ exists such that $d(T\mathcal F, T\mathcal G)<\epsilon$.  For brevity, this is stated as, $\mathcal G$ can be found {\sl $C^0$ close} to $\mathcal F$.  An isotopy $\mathcal F_t$ of $\mathcal F$ is called {\sl $C^0$ small} if it can be chosen so that at all times $\mathcal F_t$ is $C^0$ close to $\mathcal F$.  An isotopy can be found {\sl $C^0$ close to the identity} if given any $\epsilon >0$ an isotopy can be found that keeps every section within $\epsilon$ of its starting position.  

Throughout the paper, $I$ will be used to denote the closed interval $[0,1]$.  

\section{Flow boxes}

Flow box decompositions were introduced and shown to exist in \cite{KR2}.  In the definition given below, an extra condition, (5), is added that is particularly useful for inductive arguments.

\begin{definition} \cite{KR2} \label{flowboxdefn} 
Let $\mathcal F$ be either a $C^k$ or $C^{k,l}$ foliation, and let $\Phi$ be a smooth transverse flow.  A {\sl flow box}, $F$, is an $(\mathcal F,\Phi)$ compatible closed chart, possibly with corners.  That is, it is a submanifold diffeomorphic to $D\times I$, where $D$ is either a closed $C^k$ disk or polygon (a closed disk with at least three corners), $\Phi$ intersects $F$ in the arcs $\{(x,y)\}\times I$, and each component of $D\times \partial I$ is embedded in a leaf of $\mathcal F$.  The components of $\mathcal F\cap F$ give a family of $C^k$ graphs over $D$.  

In the case that $D$ is a polygon, it is often useful to view the disk $D$ as a 2-cell with $\partial D$ the cell complex obtained by letting the vertices correspond exactly to the corners of $D$.  Similarly, it is useful to view the flow box $F$ as a 3-cell possessing the product cell complex structure of $D\times I$.  Then $\partial_h F$ is a union of two (horizontal) 2-cells and $\partial_v F$ is a union of $c$ (vertical) 2-cells, where $c$ is the number of corners of $D$.  In the case that $D$ has no corners, we abuse language slightly and consider $\partial_h F$ to be a union of two (horizontal) 2-cells and $\partial_v F$ to be a single vertical face, where the face is the entire vertical annulus $\partial D\times I$.

Suppose $V$ is either empty or else a compact, $(\mathcal F,\Phi)$ compatible, codimension-0 submanifold of $M$.  A {\sl flow box decomposition of $M$ rel $V$}, or simply {\sl flow box decomposition of $M$}, if $V=\emptyset$, is a decomposition of $M\setminus \text{int}(V)$ as a finite union $M = V\cup (\cup_{i=1}^nF _i)$ where

\begin{enumerate}
\item each $F _i$ is a flow box,

\item $V\cap F_i$ is a union, possibly empty, of horizontal subsurfaces and vertical 2-cells of $F_i$, together possibly with some 0- and 1-cells,

\item the interiors of $F _i$ and $F _j$ are disjoint if $i \neq j$, 

\item if $i \ne j$ and $F_i \cap F_j $ is nonempty, it must be homeomorphic to a point, an interval, or a disk that is wholly contained either in $\partial_h F_i \cap \partial_h F_j$ or in a single face in each of $\partial_v F_i$ and $\partial_v F_j$, and

\item if $\Delta$ is a vertical 2-cell of $F_{n}$ and the interior of $\Delta$ intersects a vertical 2-cell $\Delta'$ of some $F_i$ with $i<n$, then $\Delta \subset \Delta'$.  
\end{enumerate}
\end{definition}

Most of the results proved in this paper use flow box decompositions relative to an empty codimension-0 submanifold.  The general definition is particularly useful for approximating foliations by contact structures, as described in \cite{KR2,KR3}, and it appears in support of that work in Corollary~\ref{Theorem5.2}.

\begin{prop}\label{existence of fbd} Suppose $\mathcal F$ is either a $C^k$ or a $C^{k,l}$ foliation of a compact manifold $M$, and let $\Phi$ be a smooth flow transverse to $\mathcal F$.  Suppose $V$ is either empty or else a compact, $(\mathcal F,\Phi)$ compatible, codimension-0 submanifold of $M$.  Then $M$ has a flow box decomposition rel $V$.  Moreover, any flow box decomposition of $V$ can be extended to a flow box decomposition of $M$.
\end{prop}

\begin{proof} Conditions (1)--(4) follow from Proposition~4.4 of \cite{KR2}.  Thus it is enough to show that a flow box decomposition satisfying (1)--(4) can be inductively subdivided so that (5) is satisfied.  

To do this, consider the union, $X$, of all vertical 2-cells contained in an $F_i$ with $i<n$ that intersect the interior of some vertical 2-cell of $F_n$.  Split $F_n$ along a finite collection of leaves of $\mathcal F \cap F_n$ that contain $(\partial_h X) \cap F_n$, and let $F_n^j$ be the resulting components.  Redefine the polygonal structure on each $F_n^j$ by decreeing that, in addition to the original vertical edges, every component of $\partial_v X \cap F_n^j$ is also a vertical edge.  

Replacing $F_n$ by the $F_n^j$ completes the inductive step of the construction.
\end{proof}

A flow box decomposition is called {\sl $V$-transitive}, or {\sl transitive}, when $V=\emptyset$, if $V_0=V$, $V_i = V_{i-1} \cup F_i$, and for $i=1,\dots, n$,
\begin{itemize}

\item[(6)] $V_{i-1}\cap F_i$ contains a vertical 2-cell of $F_i$.
\end{itemize}

Condition (6) is used in \cite{KR2, KR3} where flow boxes were needed to laterally propagate an approximating contact structure from $V$ to the rest of $M$.

\begin{prop}[Proposition~4.4, \cite{KR2}]\label{transitiveflowbox} If $M$ is compact and each point in $M$ can be reached from $V$ by a path in a leaf of $\mathcal F$, then there is a transitive flow box decomposition of $M$ rel $V$.  \qed
\end{prop}

If $\mathcal B =\mathcal B(\mathcal F,\Phi)$ is a flow box decomposition of $M$ rel $V$, an isotopy of $M$ is {\sl $\mathcal B$ compatible} if it is $\Phi$ compatible and, in addition, maps each cell of each flow box of $\mathcal B$ to itself setwise.

By Condition~(5), the set of vertical faces of the flow boxes $F_i$ is partially ordered by set containment; if $\Delta_1$ and $\Delta_2$ are vertical faces of $F_i$ and $F_j$, respectively, and their interiors have nonempty intersection, then $\Delta_1\subseteq \Delta_2$ or $\Delta_2\subseteq \Delta_1$.  Call a vertical face {\sl maximal} if it is maximal with respect to this partial ordering; namely, if it is not properly contained in any vertical face.

Let $\sigma_1,...,\sigma_m$ be a listing of the maximal faces.  It will sometimes be helpful to consider a regular neighborhood of $\cup_j \sigma_j$ of the following sort.

\begin{definition} \label{rnd defn} Let $F_1,...,F_n$ be a listing of the flow boxes of a flow box decomposition $\mathcal B=\mathcal B(\mathcal F,\Phi)$.  A {\sl regular neighborhood structure} $\mathcal N_{\mathcal B}=\mathcal N_{\mathcal B}(\mathcal F,\Phi)$ for $\mathcal B$ is a tuple of the form
$$(N, N_v, N(\sigma_1),...,N(\sigma_m)),$$
where 
\begin{enumerate}

\item $\sigma_1,...,\sigma_m$ is a listing of the maximal faces of $\mathcal B$,

\item each $N(\sigma_j)$ is a flow box that properly contains $\sigma_j$,

\item $N_v$ is a $(\mathcal F,\Phi)$ compatible regular neighborhood of the union of the vertical 1-cells of the maximal faces $\sigma_j$,

\item $N_v$ decomposes as a finite union of flow boxes $B_p=D_p\times I$, where $B_p\cap B_q\subset \partial_h B_p\cap \partial_h B_q$, and $B_p\cap (\cup_i (\partial_v F_i)^{(1)})\subset (\{0\}\times I)$ for each $p$,

\item if $j\ne k$, then $N(\sigma_j)\cap N(\sigma_k)$ is contained in the interior of $N_v$,

\item $N=N_v\cup_j N(\sigma_j)$, and

\item $\cup_i \partial_v F_i$ is a deformation retract of $N$.  

\end{enumerate}
\end{definition}

\begin{figure}[htbp] 
 \centering
 \includegraphics[width=2.7in]{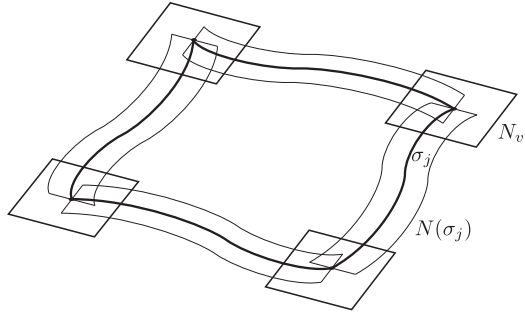} 
 \caption{Horizontal cross-section of a flow box}
 \label{regnbdsection}
\end{figure}

Figure~\ref{regnbdsection} illustrates a horizontal cross section of a regular neighborhood structure in a neighborhood of a single flow box.  Note that a $\mathcal B$ compatible isotopy takes a regular neighborhood structure $\mathcal N_{\mathcal B}(\mathcal F,\Phi)$ to a regular neighborhood structure $\mathcal N_{\mathcal B}(\mathcal F',\Phi)$.

A standard method of proof is to work inductively with a cell complex, working first with 0-cells, and then extending over the 1-cells, followed by the 2-cells, and finally the 3-cells.  When smoothness is a priority, it is often useful to work instead with regular neighborhoods of the cells.  Regular neighborhood structures provide a vocabulary for this approach in the context of flow box decompositions; namely, establish a property first on $N_v$, then on the union $\cup_j N(\sigma_j)$, and finally extend this property over the 3-cells complementary to $N\cup_i \partial_h F_i$.  

\begin{definition} A flow box decomposition $\mathcal B$ is {\sl smooth sided} if the interior of every vertical face of every flow box $F_i$ of $\mathcal B$ is a smooth surface.  The flow box decomposition is called {\sl smooth} if it is smooth-sided and every horizontal face has a neighborhood in the leaf it is contained in that is smoothly embedded.
\end{definition}

\begin{lemma}\label{smooth-sided} Let $M$ be compact.  If $\mathcal F$ is $C^{1,0}$ and $\Phi$ is a smooth transverse flow in $M$, then there exists a smooth-sided flow box decomposition of $M$.  If $\mathcal F$ is $C^{\infty,0}$ and $\Phi$ is a smooth transverse flow in $M$, then there exists a smooth flow box decomposition of $M$.
\end{lemma}

\begin{proof} The proof of Proposition~4.4 of \cite{KR2} starts by choosing initial flow boxes and these may be taken to be smooth sided.  The rest of the construction involves transversality of vertical intersections and splitting along leaves, and both of these operations work with smooth vertical faces.
\end{proof}

\begin{lemma}\label{smooth fbd} If $\mathcal F$ is $C^{1,0}$, $\Phi$ is a smooth transverse flow in $M$, and $\mathcal B$ is a smooth-sided flow box decomposition, then there exists a flow compatible isotopy that takes $\mathcal F$ to a $C^0$ close $C^{1,0}$ foliation and takes $\mathcal B$ to a smooth flow box decomposition.
\end{lemma}

\begin{proof} Let $U$ be the union of small neighborhoods, in leaves of $\mathcal F$, of each of the horizontal faces of all $F_i \in \mathcal B$.  Then $U$ is a $C^1$ embedded surface.  This may be isotoped, while preserving flow lines of $\Phi$ and keeping $\mathcal F$ $C^0$ $\epsilon$ close to itself, to a smoothly embedded surface.  Applying this isotopy to $\mathcal B$ produces the desired smooth flow box decomposition.
\end{proof}

\begin{remark} Suppose $\mathcal B$ is a smooth $(\mathcal F,\Phi)$ flow box decomposition, where $\mathcal F$ is a $C^{k,0}$ foliation for some $k\ge 1$, and $\Phi$ is a smooth flow transverse to $\Phi$.  Each flow box describes a smooth chart for $M$ in which the flow restricts to the union of vertical segments $\{\overline{x}\}\times I$ and the leaves of $\mathcal F$ restrict to a $C^0$ family of $C^k$ graphs.  After fixing a point $\bar x_0$, an index $t$ may be chosen so that the leaf containing $(\bar x_0,t)$ is given by the graph $z=f_t(\overline{x})$.  
\end{remark}

We now give several elementary and frequently used smoothing operations that will be used in a neighborhood of a surface.  To streamline statements, let $S$ be a surface, possibly with boundary, and let $S \times I \subset M$.  A {\sl strictly horizontal foliation of $S\times I$} is the foliation with leaves $S\times \{t\}, t\in I$.  An {\sl almost horizontal foliation} of $S\times I$ is a foliation transverse to the $I$ fibers which contains $S \times \partial I$ as leaves.  A product submanifold, $S\times I$ of $M$, is called an {\sl $(\mathcal F,\Phi)$ compatible product} if the restriction of $\mathcal F$ to $S\times I$ is almost horizontal, and the $I$-fibers $\{x\}\times I$ are flow segments of $\Phi$.

\begin{notation} \label{nbd notation} If $S$ is a proper subsurface of a leaf of a given foliation, $N(S)$ will denote a regular neighborhood of $S$ in its leaf.
\end{notation}

Denote by $N(S) \times I$ an extension of the product structure on $S \times I$.  If $S\times I$ is a $(\mathcal F,\Phi)$ compatible product, $N(S)\times I$ denotes a $(\mathcal F,\Phi)$ compatible product.  For some estimates, it is necessary to consider a smooth damping function $\ell:N(S) \to I$ that vanishes on $\partial N(S)$ and is identically $1$ on $S$.  

The results established in the rest of this section are elementary, useful, and substantially similar.  

They can be grouped as existence of approximations, existence of extensions, or uniqueness results, and they all occur in the context of $C^0$ approximation.  These results are typically applied when $S$ is a surface and $S\times I$ is a subset of $M$.  Given a $C^{k,r}$ foliation on all, or just a portion of $S \times I$, we ask if it can be extended to all of $S \times I$, if it can be approximated by a smoother foliation, and to what extent the approximating foliation is unique.

Complicating the statements, though not the proofs, many of the results when applied require a relative version.  Vertically, when $S$ has a boundary, the goal is to not change the given foliation near $\partial S \times I$.  Horizontally, there may be leaves on which a given foliation is as smooth as needed and should not be changed.  

In the next proposition, we establish that $C^{1,0}$ almost horizontal foliations of $S\times I$ are $C^0$ close to $C^{\infty}$ almost horizontal foliations.
 
\begin{prop} [Smoothing a product foliation] \label{smooth interpolation} Let $S$ be a compact smooth surface, and let $\mathcal P$ be a $C^{1,0}$ almost horizontal product foliation on $S \times I$.  Then $\mathcal P$ can be $C^0$ deformed to a $C^0$ close, smooth, almost horizontal product foliation $\mathcal S$ on $S \times I$.  If $\mathcal P$ is smooth on some compact $(\mathcal P,\Phi)$ compatible submanifold, then we may choose the deformation to fix this submanifold pointwise.  In addition, if some finite number of leaves of $\mathcal F$ are smoothly embedded, then the deformation can be chosen to fix these leaves pointwise, and if $L'_1,...,L'_n$ are subsurfaces of leaves $L_1,...,L_n$ of $\mathcal F$ so that regular neighborhoods $N(L'_i)$ of $L'_i$ in $L_i$ are smoothly embedded in $M$, then the deformation can be chosen to fix each $L_i$ pointwise.
\end{prop}

\begin{proof} Pick a metric on the bundle of tangent two planes to $S \times I$.  Fix a point $x_0$ in $S$, and denote the leaf of $\mathcal P$ that contains $(x_0,t)$ by $P_t$.  Given $s<t$ and $x\in S$, let $[s,t]_x$ denote the subinterval of $\{x\} \times I$ with boundary points in $P_s \cup P_t$.  Each $P_t$ is the graph of a $C^1$ function $f_t:S \to I$.  It is enough to deform the continuously varying family $f_t$ to a smoothly varying family whose graphs foliate $S \times I$.

Let $\epsilon > 0$.  Choose a partition $0=t_0< t_1<...<t_n=1$ of $[0,1]$ with the property that for all $x\in S$, for each i, the tangent planes to leaves of $\mathcal P$ at each point of $[t_{i-1},t_i]_x$ are all within $\epsilon$ of each other.  Perform a $C^0$ small isotopy so that $\cup_iP_{t_i}$ is smoothly embedded.  If some finite number of leaves of $P$ are smoothly embedded or contain subsurfaces smoothly embedded in $M$, the partition can be chosen so that these leaves appear as $P_{t_i}$.

Let $\ell_i:[t_{i-1},t_i] \to [0,1]$ be a smooth bijection that vanishes to infinite order at the endpoints.  For $t \in [t_{i-1},t_i]$ define $$g_t = (1-\ell_i(t))f_{t_{i-1}} + \ell_i(t)f_{t_i}.$$
Then $g_t$ is a smooth family of functions whose graphs give a smooth foliation $\mathcal S$ of $S\times I$.  Since $g_t$ is a linear combination of $f_{t_{i-1}}$ and $f_{t_i}$, it is easily checked using local coordinates on $S$ that the normal vector to the graph of $g_t$ is a linear combination of normals to $f_{t_{i-1}}$ and $f_{t_i}$.  It follows that tangent planes to $g_t$ $C^0$ approximate the tangent planes to $f_t$.

Since $f_t$ and $g_t$ are graphs, there is an $I$ fiber preserving deformation of $\mathcal P$ to $\mathcal S$.
\end{proof}

Sometimes, the ``smoothing'' of a horizontal foliation is required to preserve an existing structure, and $C^{\infty,0}$ smoothing is the best that can be hoped for.  The proof of Proposition~\ref{smooth interpolation} modifies easily to yield the following.

\begin{prop} [$C^{\infty,0}$ smoothing a product foliation] \label{smoothleaf interpolation} Let $S$ be a compact smooth surface, and let $\mathcal P$ be a $C^{1,0}$ almost horizontal product foliation on $S \times I$.  Then $\mathcal P$ can be $C^0$ deformed to a $C^0$ close, $C^{\infty,0}$, almost horizontal product foliation $\mathcal S$ on $S \times I$.  If $\mathcal P$ is $C^{\infty,0}$ on some compact $(\mathcal P,\Phi)$ compatible submanifold, then we may choose the deformation to fix this submanifold pointwise.  In addition, if some finite number of leaves of $\mathcal F$ are smoothly embedded, then the deformation can be chosen to fix these leaves pointwise, and if $L'_1,...,L'_n$ are subsurfaces of leaves $L_1,...,L_n$ of $\mathcal F$ so that regular neighborhoods $N(L'_i)$ of $L'_i$ in $L_i$ are smoothly embedded in $M$, then the deformation can be chosen to fix each $L_i$ pointwise.
\qed
\end{prop}

\begin{cor} [Local product smoothing in $M$]\label{local smoothing} Let $\mathcal F$ be a $C^{1,0}$ (respectively, $C^{\infty,0}$) foliation of $M$, and $\Phi$ a smooth flow transverse to $\mathcal F$.  Suppose that $S\times I$ is a $(\mathcal F,\Phi)$ compatible product, smoothly embedded in $M$, and fix $N(S)$.  Then there is a $C^0$ small, $\Phi$-compatible, isotopy of $M$ that is the identity outside $N(S)\times I$, and takes $\mathcal F$ to a $C^0$ close $C^{1,0}$ (respectively, $C^{\infty,0}$) foliation that is $C^{\infty}$ on $S\times (0,1)$.  If $\mathcal B$ is a smooth flow box decomposition of $M$, then this isotopy can be chosen to be $\mathcal B$ compatible.  If some finite number of leaves of $\mathcal F$ are smoothly embedded, then the deformation can be chosen to fix these leaves pointwise.  Moreover, if $L'_1,...,L'_n$ are subsurfaces of leaves $L_1,...,L_n$ of $\mathcal F$ so that regular neighborhoods $N(L'_i)$ of $L'_i$ in $L_i$ are smoothly embedded in $M$, then the deformation can be chosen to fix each $L'_i$ pointwise.
\end{cor}

\begin{proof} If $M=S\times I$, this follows immediately from Proposition~\ref{smooth interpolation}.  Restrict attention, therefore, to the case that $S$ is a proper subsurface of a leaf of $\mathcal F$.

Fix a Riemannian metric on $M$.  Choose a regular neighborhood $N(S)$ of $S$, and product structure $N(S)\times I$ extending the product structure $S\times I$, so that the restriction of $\mathcal F$ to $N(S)\times I$ remains almost horizontal.  The restriction of $\mathcal F$ to $N(S)\times I$ is therefore a product foliation, $\mathcal P$ say.  Let $\ell:N(S) \to I$ be a fixed smooth damping function that vanishes on $\partial N(S)$ and is identically $1$ on $S$.

By Proposition~\ref{smooth interpolation}, $\mathcal P$ can be $C^0$ deformed to a $C^0$ close, $C^{\infty}$, almost horizontal product foliation $\mathcal S$ on $N(S) \times I$.  If some finite number of leaves of $\mathcal F$ are smoothly embedded, then the deformation can be chosen to fix these leaves pointwise.  Moreover, if $L'_1,...,L'_n$ are subsurfaces of leaves $L_1,...,L_n$ of $\mathcal F$ so that regular neighborhoods $N(L'_i)$ of $L'_i$ in $L_i$ are smoothly embedded in $M$, then the deformation can be chosen to fix each $L'_i$ pointwise.

Using the notation found in the proof of Proposition~\ref{smooth interpolation}, let $f_t$ and $g_t$, $0\le t\le 1$, be functions over $N(S)$ whose graphs give the leaves of $\mathcal P$ and $\mathcal S$, respectively.  If $f_t$ describes the graph of the leaf $L_i$ for some $i$, then $f_t=g_t$ along $L'_i$.  Note that given any $\epsilon>0$, the proof of Proposition~\ref{smooth interpolation} guarantees that we may choose the smoothly varying family $g_t$ so that both the graphs $g_t$ and $f_t$ and their tangent plane fields are $\epsilon$ $C^0$ close, for each $t\in [0,1]$.

Let $s\in [0,1]$ and define $h_t^s:N(S) \to I$ by $$h_t^s(x) = (1 - s\ell(x))f_t(x) + s\ell(x)g_t(x).$$ For each $s$, points of $N(S) \times I$ are uniquely expressible as $(x,h_t^s(x))$, thus $h_t^s$ defines a fiber preserving isotopy of $N(S) \times I$ that is the identity in a neighborhood of $\partial N(S)\times I$ and takes $\mathcal P$ to $\mathcal S$ on $S\times I$.  

Since we can choose $\mathcal S$ to be arbitrarily $C^0$ close to $\mathcal P$, we can guarantee that this isotopy is $C^0$ small, and the resulting foliation is $C^0$ close to $\mathcal F$.  To see this, let $s\in [0,1]$ and consider $$h_t^s(x) -f_t(x)=s\ell(x)(g_t(x) - f_t(x)).$$ Suppose $x$ is given locally on $S$ by real coordinates $(u,v)$.  Taking the partial with respect to $u$ gives
$$s\partial_u\ell(x)(g_t(x) - f_t(x))+ s\ell(x)\partial_u(g_t(x)-f_t(x)).$$ 
A symmetric equation gives the derivative with respect to $v$.  We start with a fixed $\ell$, and, for any $\epsilon>0$, we may choose the smoothly varying family $g_t$ so that both the graphs $g_t$ and $f_t$ and their tangent plane fields are $\epsilon$ $C^0$ close, for each $t\in [0,1]$.  Therefore, we can guarantee that both the graphs $h_t^s$ and $f_t$ and their tangent plane fields are $\epsilon$ $C^0$ close, for each $s,t\in [0,1]$.
\end{proof}

\begin{cor} [Local smoothing - uniqueness] \label{local smoothing unique} Let $\mathcal F$ be a $C^{1,0}$ foliation of $M$, and $\Phi$ a smooth flow transverse to $\mathcal F$.  Suppose that $S\times I$ is a $(\mathcal F,\Phi)$ compatible product, smoothly embedded in $M$, and fix $N(S)$.  Suppose $\mathcal S$ is an almost horizontal product foliation of $S \times I$.  There is an isotopy of $\mathcal F$, fixing the complement of $N(S) \times I$, that takes $\mathcal P$ to $\mathcal S$.  The distance that a tangent plane to $\mathcal F$ moves during this isotopy depends on the distance between the tangent planes to $\mathcal F$ and $\mathcal S$ on $N(S) \times I$ and the maximum of the derivative of the damping function used.  
\end{cor}

\begin{proof} This will follow from a more general statement, Proposition~\ref{uniqueness in a product} that allows non-trivial holonomy.
\end{proof}

\begin{definition} Let $\mathcal F$ be an almost horizontal foliation on $S \times I$.  If $\beta \subset S$ is an arc with $\beta(0)=*$ and $\beta(1)=x$, let $\rho_\mathcal F(\beta)$ denote the homeomorphism from $\{ x\} \times I \to \{*\} \times I$ given by lifting $\beta$ to leaves of $\mathcal F$.  More precisely, given such an arc $\beta$, let $\beta_t$ be the path in a leaf of $\mathcal F$ that ends at $(x,t)$ and projects to $\beta$, and define $\rho_\mathcal F(\beta)(x,t)=\beta_t(0)$.  Let $[\beta]$ be the path homotopy class of $\beta$, and define $\rho_\mathcal F([\beta])=\rho_\mathcal F(\beta)$.

Restricting attention to loops in $S$ based at $*$, we obtain the {\sl holonomy representation} $\rho_\mathcal F: \pi_1(S, *) \to \text{Homeo}(\{*\} \times I)$.  
\end{definition}

\begin{prop} [Uniqueness in a product]\label{uniqueness in a product} Let $\mathcal F$ and $\mathcal G$ be foliations of $M$ that restrict to almost horizontal foliations of $N(S) \times I$.  If the holonomy representations of $\mathcal F$ and $\mathcal G$ agree on $S \times I$, then there is a $C^0$ isotopy of $M$ that preserves $I$ bundle fibers, is supported in a neighborhood on $N(S) \times I$, and on $S \times I$ takes $\mathcal F$ to $\mathcal G$.  The distance that a tangent plane to $\mathcal F$ moves during this isotopy depends on the distance between the tangent planes to $\mathcal F$ and $\mathcal G$ on $N(S) \times I$ and the maximum of the derivative of the damping function used.
\end{prop}

\begin{proof} Fix $*$ in $S$.  For each $x\in S$, choose a path $\beta$ in $S$ from $*$ to $x$, and define a homeomorphism $h$ of $\{x\} \times I$ by $h(t) =\rho^{-1}_{\mathcal G}(\beta) \rho_{\mathcal F}(\beta)(x,t)$.  Then $h$ gives a well defined $I$ fiber preserving homeomorphism of $S\times I$ taking $\mathcal F$ to $\mathcal G$ that can be used to produce the desired straight line isotopy.
\end{proof}

In the next two propositions, we focus attention on flow boxes of the form $F=D\times I$, where $D=I\times I$, and establish leafwise smoothing subject to a monodromy constraint.  Let $\alpha:[0,1]\to D$ be the arc given by $\alpha(t)=(1/2,t)$.  Given an almost horizontal foliation $\mathcal F$ of $F$, we have 
$$\rho_{\mathcal F}(\alpha): \{(1/2,1)\} \times I \to \{(1/2,0)\} \times I.$$

\begin{prop} [$C^{\infty,0}$ smoothing a product foliation with holonomy constraint] \label{leafsmoothingholconstraints}
Let $D=I\times I$, with core curve $\alpha=\{1/2\}\times I$, oriented as above.  Let $\Phi$ denote the 1-dimensional foliation by $I$ fibers $\{x\}\times I, x\in D$.  Suppose $\mathcal P$ is a $C^{1,0}$ almost horizontal foliation on $D\times I$ that is $C^{\infty,0}$ on the neighborhood $N_v=I\times (J_0\cup J_1)\times I$ of $I\times \partial I\times I$, where $J_0$ and $J_1$ are nondegenerate closed intervals containing $0$ and $1$ respectively.  Then $\mathcal P$ can be $C^0$ deformed to a $C^0$ close, $C^{\infty,0}$, almost horizontal foliation $\mathcal G$ on $D \times I$ such that $\mathcal G$ agrees with $\mathcal P$ on $N_v$ and $\rho_{\mathcal G}(\alpha)=\rho_{\mathcal P}(\alpha)$.  If $L_1',...,L_n'$ are smoothly embedded leaves of $\mathcal P$, then the deformation can be chosen to fix these leaves pointwise.  Moreover, if $L'_1,...,L'_n$ are subsurfaces of leaves $L_1,...,L_n$ of $\mathcal P$ so that regular neighborhoods $N(L'_i)$ of $L'_i$ in $L_i$ are smoothly embedded in $D\times I$, then the deformation can be chosen to fix each $L'_i$ pointwise.
\end{prop}

\begin{proof} Let $N(N_v)=I\times N(J_0\cup J_1)\times I$ be an open neighborhood of $N_v$ on which $\mathcal P$ is $C^{\infty,0}$.

Pick a metric on the bundle of tangent 2-planes to $D\times I$.  Denote the leaf of $\mathcal P$ that contains $(1/2,0,t)$ by $P_t$.  Each $P_t$ is the graph of a $C^1$ function $f_t:D\to I$.  Fix a smooth damping function $\ell:I \to I$ that vanishes on $J_0\cup J_1$ and is identically $1$ on the complement of $ N(J_0\cup J_1)$.  

By Proposition~\ref{smoothleaf interpolation}, there is a $C^{\infty,0}$, almost horizontal, foliation $\mathcal S$ of $D\times I$ that is $C^0$ close to $\mathcal P$ and agrees with $\mathcal P$ on $N_v$.  Denote the leaf of $\mathcal S$ that contains $(1/2,0,t)$ by $S_t$.  Each $S_t$ is the graph of a $C^{\infty}$ function $s_t:D\to I$.  If $f_t$ describes the graph of the leaf $L'_i$ for some $i$, then choose $\mathcal S$ so that $s_t=f_t$ along $L'_i$.

It may be helpful, even though not necessary, to recall that by the construction of $\mathcal S$ found in the proof of Proposition~\ref{smooth interpolation}, we may assume that that there is a partition $0=t_0<t_1<...<t_m=1$, such that for each $i$, $S_{t_i}=P_{t_i}$, and each leaf $S_t, t\in [t_i,t_{i+1}],$ has tangent plane field $C^0$ close to the tangent plane field of each of $P_{t_i}$ and $P_{t_{i+1}}$.  

Define a homeomorphism $h:I\to I$ by
$$ ((1/2,0),h(t))= \rho_{\mathcal S}(\alpha)\circ \rho^{-1}_{\mathcal P}(\alpha)((1/2,0),t).$$

The point of the definition of $h(t)$ is that a leaf of $\mathcal S$ that contains $((1/2,0),h(t))$ will intersect the leaf of $\mathcal P$ that contains $((1/2,0),t)$ at the point $((1/2,1),s_{h(t)}(1/2,1))=((1/2,1),f_t(1/2,1))$.

Since $\mathcal P$ and $\mathcal S$ agree on $I \times (J_0 \cup J_1)$, $s_{h(t)} = f_t$ on all of $I \times (J_0 \cup J_1)$.  Since $\mathcal P$ and $\mathcal S$ are $C^0$ close, $h$ is $C^0$ close to the identity map, and, for each $t\in I$, the graphs $s_t$ and $s_{h(t)}$ are $C^1$ close on the rest of $D$.

We obtain a $C^{\infty,0}$ foliation $\mathcal G$ approximating $\mathcal P$ and satisfying the holonomy constraint $\rho_{\mathcal G} (\alpha)=\rho_{\mathcal P} (\alpha)$ as follows.  Define
\begin{align*}g_t(x,y) &= \ell(y) s_t(x,y)+ (1-\ell(y)) s_{h(t)}(x,y)\\
& =\ell(y)(s_t(x,y) - s_{h(t)}(x,y)) + s_{h(t)}(x,y).  
\end{align*}
 Then for a fixed $t$, 
 $$\partial_x g_t = \ell\cdot (\partial_x s_t - \partial_x s_{h(t)}) + \partial_x s_{h(t)} $$
and $$\partial_y g_t = \partial_y \ell (s_t-s_{h(t)}) + \ell\cdot (\partial_y s_t - \partial_y s_{h(t)}) + \partial_y s_{h(t)}.$$ 
By choosing $\mathcal S$ sufficiently close to $\mathcal P$, we guarantee that the $C^{\infty,0}$ foliation $\mathcal G$ with leaves the graphs of the functions $g_t:D\to I$ is $C^0$ close to $\mathcal P$.  Since $\mathcal S=\mathcal P$ on $N_v$, the choice of $\ell$ implies $\mathcal G=\mathcal P$ on $N_v$.  And, as has now been argued many times, this foliation $\mathcal G$ is obtained by a $C^0$ small deformation of $\mathcal P$.
\end{proof}

A similar argument yields the following.

\begin{prop} [$C^{\infty}$ smoothing a product foliation with smooth holonomy constraint] \label{smoothingholconstraints}
Let $D=I\times I$, with core curve $\alpha=\{1/2\}\times I$, oriented as above.  Let $\Phi$ denote the 1-dimensional foliation by $I$ fibers $\{x\}\times I, x\in D$.  Suppose $\mathcal P$ is a $C^{1,0}$ almost horizontal foliation on $D\times I$ that is smooth on the neighborhood $N_v=I\times (J_0\cup J_1)\times I$ of $I\times \partial I\times I$, where $J_0$ and $J_1$ are nondegenerate closed intervals containing $0$ and $1$ respectively.  Suppose that the holonomy map $\rho_{\mathcal P} (\alpha)$ is smooth.  Then $\mathcal P$ can be $C^0$ deformed to a $C^0$ close, smooth, almost horizontal foliation $\mathcal G$ on $D \times I$ such that $\rho_{\mathcal G}(\alpha)=\rho_{\mathcal P}(\alpha)$.  Moreover, $\mathcal G$ can be chosen to agree with $\mathcal P$ on $N_v$.  If leaves $L_1',...,L_n'$ smoothly embedded leaves of $\mathcal P$, then the deformation can be chosen to fix these leaves pointwise.  Moreover, if $L'_1,...,L'_n$ are subsurfaces of leaves $L_1,...,L_n$ of $\mathcal P$ so that regular neighborhoods $N(L'_i)$ of $L'_i$ in $L_i$ are smoothly embedded in $D\times I$, then the deformation can be chosen to fix each $L'_i$ pointwise.
\end{prop}

Once again, local smoothing can be realized by a flow compatible, $C^0$ small, isotopy of $M$.  This is made precise in the following corollary.

\begin{cor}[Local product $C^{\infty,0}$ smoothing with holonomy constraint in $M$]\label{leafsmoothingcor}
Let $\mathcal F$ be a $C^{1,0}$ (respectively, $C^{\infty,0}$) foliation of $M$, and $\Phi$ a smooth flow transverse to $\mathcal F$.  Suppose that $F=D\times I$, $D=I\times I$, is a smooth flow box in $M$, and that $\mathcal F$ is $C^{\infty,0}$ in a neighborhood $N_v$ of $I\times\partial I\times I$ .  Then there is a $C^0$ small, $\Phi$-compatible, isotopy of $M$ that is the identity outside a smooth flow box $N(D) \times I$, and takes $\mathcal F$ to a $C^0$ close $C^{1,0}$ foliation $\mathcal G$ that is $C^{\infty,0}$ on $D\times I$, agrees with $\mathcal F$ on $N_v$, and satisfies $\rho_{\mathcal F }(\alpha)=\rho_{\mathcal G }(\alpha)$.  If $\mathcal B$ is a smooth flow box decomposition of $M$, then this isotopy can be chosen to be $\mathcal B$ compatible.  \qed
\end{cor}

For the rest of this section, given a flow box $F= D \times I$, $A$ will denote a smooth regular collar neighborhood of $\partial D$ in $D$, and $N(A)$ will denote a regular neighborhood of $A$ in $D$.

The next proposition shows that a smooth annular foliation on $A\times I$ that approximates $\mathcal F$ can be extended to a smooth foliation on $F$ that approximates $\mathcal F$.

\begin{prop} [$C^0$ close damped coning]\label{damped coning} Let $\mathcal F$ be a $C^{1,0}$ foliation of $M$, and $\Phi$ a smooth flow transverse to $\mathcal F$.  Let $F=D\times I$ be a $C^{\infty}$ flow box, and suppose that $\mathcal A$ is a smooth almost horizontal product foliation of $N(A)\times I$ that is $C^0$ close to the restriction of $\mathcal F$ to $N(A)\times I$.  Then there is a smooth foliation $\mathcal G$ of $F$ that is $C^0$ close to $\mathcal F$ and restricts to $\mathcal A$ on $A\times I$.
\end{prop}

\begin{proof}
By Corollary~\ref{local smoothing}, there is a smooth foliation $\mathcal D$ of $F$ which is $C^0$ close to $\mathcal F$.  Then $\mathcal A$, $\mathcal F$, and $\mathcal D$ are all $C^0$ close on $N(A) \times I$.  Since $\mathcal A$ and $\mathcal D$ are both product foliations on $N(A) \times I$, the result follows by applying uniqueness of local smoothing, Corollary~\ref{local smoothing unique}, to $A$ and $N(A)$.  
\end{proof}

This result can be strengthened so that the approximating foliation is a smooth extension of $\mathcal F$ near the vertical boundary.  

\begin{cor} [Damped coning] \label{extendsmoothoverthreecell}
Let $\mathcal F$ be a $C^{1,0}$ foliation of $M$, and $\Phi$ a smooth flow transverse to $\mathcal F$.  
Let $F=D\times I$ be a $C^{\infty}$ flow box, and suppose that the restriction of $\mathcal F$ to $N(A)\times I$ is a $C^{\infty}$ foliation.  Then there is a $C^0$ small, flow compatible, isotopy of $F$ that is the identity on $A\times I$ and takes $\mathcal F$ to a $C^0$ close $C^{\infty}$ foliation $\mathcal G$ of $F$.  If $\mathcal F$ is smooth in a neighborhood of $D\times\partial I$ in $F$, then the isotopy can be chosen so that $\mathcal G=\mathcal F$ in a slightly smaller neighborhood of $D\times\partial I$.  If $\mathcal B$ is a smooth flow box decomposition of $M$, then this isotopy can be chosen to be $\mathcal B$ compatible.
\end{cor}

\begin{proof}
Let $\mathcal A$ be the restriction of $\mathcal F$ to $N(A)\times I$.  The result then follows immediately from Proposition~\ref{damped coning}.
\end{proof}

\begin{prop} [Leafwise smooth damped coning]\label{extendleafsmoothoverthreecell}
Let $\mathcal F$ be a $C^{1,0}$ foliation of $M$, and $\Phi$ a smooth flow transverse to $\mathcal F$.  Let $F=D\times I$ be a $C^{\infty}$ flow box, and suppose that the restriction of $\mathcal F$ to $N(A)\times I$ is a $C^{\infty,0}$ foliation.  Then there is a $C^0$ small flow compatible isotopy of $F$ that is the identity on $A\times I$ and takes $\mathcal F$ to a $C^0$ close $C^{\infty,0}$ foliation $\mathcal G$ of $F$.  If $\mathcal F$ is smooth in a neighborhood of $D\times\partial I$ in $F$, then the isotopy can be chosen so that $\mathcal G=\mathcal F$ in a slightly smaller neighborhood of $D\times\partial I$.  If $\mathcal B$ is a smooth flow box decomposition of $M$, then this isotopy can be chosen to be $\mathcal B$ compatible.
\end{prop}

\begin{proof} Recall that $A$ denotes a smooth regular collar neighborhood of $\partial D$ in $D$, and $N(A)$ denotes a smooth regular neighborhood of $A$ in $D$.  Decompose $N(A)=A\cup B\cup A'$ as a union of smooth annuli with pairwise disjoint interiors and $B$ disjoint from $\partial N(A)$.  By assumption, the restriction of $\mathcal F$ to $N(A)\times I$, and hence to $A'\times I$, is a $C^{\infty,0}$ product foliation by almost horizontal annuli.

Apply Proposition~\ref{local smoothing} with $S=A'$: there is a $C^0$ small, flow compatible isotopy of $F$ that is the identity outside a small neighborhood $N(A')\times I$ and takes $\mathcal F$ to a $C^0$ close $C^{\infty,0}$ foliation $\mathcal G$ that is $C^{\infty}$ on $A'\times I$.  By choosing $N(A')$ disjoint from $A$, we may assume this isotopy is the identity on $A\times I$.  It follows that the restriction of $\mathcal G$ to $A\times I$ agrees with the restriction of $\mathcal F$ to $A\times I$.

Next, apply Proposition~\ref{damped coning} to the flow box $F'=(D\setminus \text{int}(A\cup B))\times I$.  The result is a smooth foliation $\mathcal G'$ on $F'$ that is $C^0$ close to $\mathcal F$ and agrees with $\mathcal G$ on $A'\times I$.  It therefore extends as $\mathcal G$ to give a $C^{\infty,0}$ foliation of $F$ that is $C^0$ close to $\mathcal F$ and agrees with $\mathcal F$ on $A\times I$.
\end{proof}

\section{Any $C^{1,0}$ foliation is a limit of $C^{\infty,0}$ foliations}

The next theorem adds $C^0$ approximation to a theorem of Calegari, \cite{calegari}.  It is applied and cited as Theorem~2.10 in \cite{KR3}.

\begin{thm} \label{calegarikr}
Suppose $\mathcal F$ is a $C^{1,0}$ foliation of a compact manifold $M$.  Then there is a $C^0$ small isotopy of $M$ taking $\mathcal F$ to a $C^{\infty,0}$ foliation $\mathcal G$ that is $C^0$ close to $\mathcal F$.  If $\Phi$ is a smooth flow transverse to $\mathcal F$, the isotopy may be taken to be flow compatible.
\end{thm}

\begin{proof}
Let $\Phi$ be a smooth flow transverse to $\mathcal F$, and apply Lemma~\ref{smooth-sided} to obtain a smooth-sided $(\mathcal F,\Phi)$-flow box decomposition $\mathcal B'$ of $M$.  Let $\mathcal B$ denote the smooth flow box decomposition and $\mathcal F_1$ the foliation that results from applying the isotopy of Lemma~\ref{smooth fbd} to $\mathcal B'$ and $\mathcal F$ respectively.  Let $\sigma_1,...,\sigma_n$ be a listing of the maximal vertical faces of $\mathcal B$, and choose a regular neighborhood structure $(N,N_v,N(\sigma_1),...,N(\sigma_n))$ for $\mathcal B$.

By Corollary~\ref{local smoothing}, there is a $C^0$ small, $\mathcal B$ compatible isotopy of $M$ that takes $\mathcal F_1$ to a $C^0$ close $C^{1,0}$ foliation $\mathcal F_2$ that is smooth on $N_v$.  

By Corollary~\ref{leafsmoothingcor}, there is a $C^0$ small, $\mathcal B$ compatible isotopy of $M$ that takes $\mathcal F_2$ to a $C^0$ close $C^{1,0}$ foliation $\mathcal F_3$ that is smooth on $N_v$ and $C^{\infty,0}$ on each $N(\sigma_i)$.  

Finally, by Proposition~\ref{extendleafsmoothoverthreecell}, there is a $C^0$ small, $\mathcal B$ compatible isotopy of $M$ that takes $\mathcal F_3$ to a $C^0$ close $C^{\infty,0}$ foliation $\mathcal G$.
\end{proof}

\begin{cor}\label{limitsmooth}
Any $C^{1,0}$ foliation is a limit of $C^{\infty,0}$ foliations.\qed
\end{cor}

\section{Measured Foliations}
 
A {\sl transverse measure} on a codimension one foliation $\mathcal F$ is a {\sl continuous, non-degenerate, invariant} measure, $\mu$, on each arc transverse to $\mathcal F$.  It is continuous in the sense that if $\tau$ is smoothly parametrized as $\tau=[0,x]$, then $\mu([0,x])$ is continuous in $x$.  Non-degenerate means that $\mu$ is positive on every open interval.  Invariant, in this context, means that the measure of a transverse arc is unchanged under isotopies of the arc that keep each point on the same leaf of $\mathcal F$.

\begin{lemma}[Smoothing a measure near a transversal] \label{smoothtransversal} Let $(\mathcal F,\mu)$ be a $C^{\infty,0}$ measured foliation in $M$.  Suppose $\tau$ is a smoothly embedded arc or closed curve which is everywhere transverse to $\mathcal F$.  Then there is a $C^0$ small isotopy of $M$ which is the identity outside some small regular neighborhood $N$ of $\tau$ and takes the measured foliation $(\mathcal F,\mu)$ to a $C^{\infty,0}$ measured foliation $(\mathcal F',\mu')$ such that $\mathcal F'$ is smooth in a neighborhood of $\tau$ and the measure, $\mu'$, restricted to $\tau$ is smooth.  If $\mu$ is smooth on a closed submanifold $A$ of $\tau$, then the isotopy can be chosen so that $\mu'=\mu$ on $A$.  If, in addition, $\mathcal F$ is smooth in a $(\mathcal F,\Phi)$ compatible regular neighborhood $N_0$ of $A$, then the isotopy can be chosen to be the identity on $N_0$.
\end{lemma}

\begin{proof} It suffices to consider the case that $\tau$ is a smoothly embedded arc.  

Regard $\tau$ as a smooth map $I \to M$.  Then $\mu(\tau[0,t])$ is a homeomorphism, $h:I\to\mathbb R$, onto its image.  Approximate $h$ by a diffeomorphism, relative to end points, $g$.  The goal is to make a small, continuous change of coordinates on $\tau$ so that $\mu$ is smooth in the new coordinates.  In other words, we must choose a homeomorphism $f:I \to I$ so that $h\circ f$ is smooth.  This is accomplished by defining $f = h^{-1} g$.

The next step is to use this reparametrization to describe a $C^0$ small isotopy of $M$ which is the identity outside a small neighborhood of $\tau$ and takes the measure $\mu$ to a measure $\mu'$ that on $\tau$ satisfies $\mu'[0,t] =\mu (\tau\circ f[0,t])$ for each $t\in [0,1]$.  In particular, $\mu'$ restricted to $\tau$ is smooth.  

To accomplish this, let $N_1$ and $N_2$ be small smoothly embedded tubular neighborhoods of $\tau$ satisfying $\overline{N}_1\subset \mbox{int}(N_2)$.  Choose these neighborhoods small enough so that $\mathcal F$ meets each in a foliation by meridian disks.  Parametrize these disks by their intersection with $\tau=[0,t]$ and so that $N_1$ is identified with the smooth family of meridian disks $D_t=D^2\times \{t\}$.

Next isotope $\mathcal F$ in $N_2$ so that the disk $D_{f(t)}$ is taken to the disk $D_t$.  If we choose $N_1$ small enough and $f$ sufficiently close to the identity, we may choose these isotopies to be as close as desired to the identity.

Define $\mu'$ along $\tau$ so that is invariant and agrees with $\mu$ away from $\tau$.  
\end{proof}

\begin{cor}[Smoothing $\mathcal F$ near a transversal] \label{smoothtransversal2}
Let $\mathcal F$ be a $C^{\infty,0}$ foliation in $M$.  Suppose $\tau$ is a smoothly embedded arc or closed curve which is everywhere transverse to $\mathcal F$.  Then there is a $C^0$ small isotopy of $M$ that is the identity outside some small regular neighborhood of $\tau$ and takes $\mathcal F$ to a $C^{\infty,0}$ $\mathcal F'$ such that $\mathcal F'$ is smooth in a neighborhood of $\tau$.  If $\mathcal F$ is smooth in a $(\mathcal F,\Phi)$ compatible regular neighborhood $N_0$ of a closed subset of $\tau$, then the isotopy can be chosen to be the identity on $N_0$.
\end{cor}

\begin{proof}
Choose a small regular neighborhood $N$ of $\tau$ so that $\mathcal F$ meets it in a product foliation by disks.  Use distance along $\tau$ to define a smooth transverse measure on the restriction of $\mathcal F$ to $N$.  The result now follows immediately from Lemma~\ref{smoothtransversal}.  
\end{proof}

The next lemma shows how the existence of a transverse measure allows the foliation to be smoothed near a compact portion of a leaf.  

\begin{lemma} [Smoothing in the neighborhood of a compact subsurface of a leaf] \label{smoothhorizontal} Let $(\mathcal F,\mu)$ be a $C^{\infty,0}$ measured foliation in $M$.  Suppose $S$ is a compact subsurface of a leaf of $\mathcal F$.  Then there is a $C^0$ small isotopy of $M$ which is the identity outside some small regular neighborhood of $S$ in $M$ and takes the measured foliation $(\mathcal F,\mu)$ to a $C^{\infty,0}$ measured foliation $(\mathcal F',\mu')$ such that $\mathcal F'$ is smooth in a neighborhood of $S$ and the measure $\mu'$ restricted to this neighborhood is smooth.
\end{lemma}

\begin{proof} Let $L$ be the leaf of $\mathcal F$ containing $S$.  If $S=L$, let $N(S)=L$.  Otherwise, let $N(S)$ be the closure of a regular neighborhood of $S$ in $L$.  Use the measure, $\mu$, to give a continuous parametrization of the flow $\Phi$ in a neighborhood of $ N(S)$.  To avoid confusion, let $\Phi'$ denote this reparametrized restriction of $\Phi$.  Choose this parametrization so that $\Phi'(x,0)=x$ for all $x\in N(S)$ and $\mu(\Phi'(x,[s,t]))=t-s$ for $s<t$ sufficiently close to $0$.  

For some $\epsilon>0$, $\Phi':N(S)\times [-\epsilon,\epsilon]\to M$ is a topological embedding, and for each $t\in [-\epsilon,\epsilon]$, $\Phi'(N(S)\times\{t\})$ is a compact subsurface of a leaf of $\mathcal F$, necessarily isotopic to $N(S)$.  Since $\mathcal F$ is $C^{\infty,0}$ and $\Phi'$ is smooth when restricted to a leaf, $\Phi'(N(S)\times [-\epsilon,\epsilon] )$ is a smooth codimension 0 submanifold, possibly with corners.  

Use Proposition~\ref{smooth interpolation} to $C^0$ isotope $\mathcal F$ in $\Phi'(N(S)\times [-\epsilon,\epsilon] )$ so that it is a smooth foliation by surfaces isotopic to $N(S)$.  The resulting measured foliation $(\mathcal F',\mu')$ and the measure $\mu'$ are necessarily smooth on the neighborhood $\Phi'(N(S)\times (-\epsilon,\epsilon))$ of $N(S)$.  
\end{proof}

The next theorem is applied and cited as Theorem~8.10 in \cite{KR3}:

\begin{thm}\label{msredtosmooth0} Suppose $\mathcal F$ is a transversely orientable $C^{1,0}$ measured foliation in $M$.  Then there is an isotopy of $M$ taking $\mathcal F$ to a $C^{\infty}$ measured foliation which is $C^0$ close to $\mathcal F$.  If $\Phi$ is a smooth flow transverse to $\mathcal F$, the isotopy may be taken to be flow compatible.
\end{thm}

\begin{proof} 
By Theorem~\ref{calegarikr} we may assume $\mathcal F$ is $C^{\infty,0}$.  From Lemma~\ref{smooth-sided} it follows that $M$ has a smooth flow box decomposition, $M= F_1 \cup \dots \cup F_n$.

Using Lemma~\ref{smoothhorizontal}, we may assume, after a $C^0$ small isotopy, that $\mathcal F$ and $\mu$ are smooth in a small regular neighborhood $N_h$ of $\cup_i \partial_h F_i$.  The vertical 1-skeleton of $\mathcal B$ is a disjoint union of transversals to $\mathcal F$, and hence, by applying Lemma~\ref{smoothtransversal}, we may assume, after a $C^0$ small isotopy, that $\mathcal F$ and $\mu$ are smooth on $N_h\cup N_v$, where $N_v$ is a union of flow boxes that form a $(\mathcal F,\Phi)$ compatible regular neighborhood of the union of the vertical 1-cells of $\cup_i \partial^{(1)} F_i$.  

Let $\sigma_1,...,\sigma_n$ be a listing of the maximal vertical faces of $\mathcal B$, and let $(N,N_v,N(\sigma_1)...,N(\sigma_n))$ be a regular neighborhood structure for $\mathcal B$.  Choose $N(N(\sigma_i)), 1\le i\le n$, so that $$( N \cup \cup_iN(N(\sigma_i)),N_v,N(N(\sigma_1)),...,N(N(\sigma_n)))$$ is also a regular neighborhood structure for $\mathcal B$.  

Let $\tau_i^-$ and $\tau_i^+$ denote the two vertical edges of $\sigma_i$.  Let $\rho_i$ denote the homeomorphism obtained by following leaves of $\mathcal F$ across $\sigma_i$.  Since the measure $\mu$ is smooth on $\tau_i^{\pm}$, and $\rho_i$ preserves $\mu$, $\rho_i$ is smooth.  By Proposition~\ref{smoothingholconstraints}, there is a $C^0$ small, $\mathcal B$ compatible isotopy of $M$ that is the identity on $N_v\cup N_h$ and outside the union $\cup_i N(N(\sigma_i))$, and takes $\mathcal F$ to a foliation $\mathcal G$ that is smooth on $N$ and $C^0$ close to $\mathcal F$.

Finally, we apply damped coning, Corollary~\ref{extendsmoothoverthreecell}, to extend $\mathcal G$ to a smooth foliation that is $C^0$ close and isotopic, by a $C^0$ small isotopy, to $\mathcal F$.  Since $\mu$ is defined on the vertical boundary of every flow box, it extends to all of $\mathcal G$.
\end{proof}

The next corollary now follows from Theorem~{8.11} in \cite{KR3}, which uses Theorem~\ref{msredtosmooth0} in its proof.  Alternatively, it also follows from Theorem~\ref{msredtosmooth0} together with Tischler's Theorem \cite{Ti}, which states that any transversely oriented, measured $C^2$ foliation on a compact n-manifold can be $C^{\infty}$ approximated by a smooth fibration over $S^1$.

\begin{cor} \label{Tischler} A $C^{1,0}$, transversely oriented, measured foliation on a compact 3-manifold is $C^0$ close to a smooth fibration over $S^1$.  \qed 
\end{cor}

\section{Holonomy Neighborhoods}\label{holonomysection}

In this section, we recall some definitions and results used in \cite{KR3}, giving those proofs which, for clarity of exposition, were deferred to this paper.  

Let $\gamma$ be an {\sl oriented} simple closed curve in a leaf $L$ of $\mathcal F$, and let $p$ be a point in $\gamma$.  We are interested in the behavior of $\mathcal F$ in a neighborhood of $\gamma$.  Let $h_{\gamma}$ be a holonomy map for $\mathcal F$ along $\gamma$, and let $\sigma$ and $\tau$ be small closed segments of the flow $\Phi$ which contain $p$ in their interiors and satisfy $h_{\gamma}(\tau) = \sigma $.  Choose $\tau$ small enough so that $\sigma \cup \tau$ is a closed segment and not a loop.  Notice that $\sigma\cap\tau$ is necessarily a closed segment containing $p$ in its interior.  There are three possibilities:
\begin{enumerate}
\item $\sigma=\tau$,
\item one of $\sigma$ and $\tau$ is properly contained in the other, or
\item $\sigma\cap\tau$ is properly contained in each of $\sigma$ and $\tau$.  
\end{enumerate}

We will need to consider very carefully a regular neighborhood of $\gamma$ which lies nicely with respect to both $\mathcal F$ and $\Phi$.  To this end, restrict attention to foliations $\mathcal F$ which are $C^{\infty,0}$ and transversely oriented, and transverse flows $\Phi$ which are smooth, and suppose that $\gamma$ is smoothly embedded in $L$.  Let $A$ be the closure of a smooth regular neighborhood of $\gamma$ in $L$; so $A$ is a smoothly embedded annulus in $L$.

\begin{lemma} [Lemma~3.1 of \cite{KR3}] \label{holnbddefn} Suppose $\mathcal F$ is $C^{\infty,0}$ and transversely oriented, and $\Phi$ is smooth.  If $\tau$ and $A$ are chosen to be small enough, there is a compact submanifold $V$ of $M$, smoothly embedded with corners, which satisfies the following:
\begin{enumerate}
\item $V$ is homeomorphic to a solid torus,
\item $\partial V$ is piecewise vertical and horizontal; namely, $\partial V$ decomposes as a union of subsurfaces $\partial_v V\cup \partial_h V$, where $\partial_v V$ is a union of flow segments of $\Phi$ and $\partial_h V$ is a union of two surfaces $L_-$ and $L_+$, each of which is either a disk or an annulus, contained in leaves of $\mathcal F$,
\item each flow segment of $\Phi\cap V$ runs from $L_-$ to $L_+$,
\item $\tau$ is a component of the flow segments of $\Phi\cap V$, and
\item $A$ is a leaf of the foliation $\mathcal F\cap V$.
\end{enumerate}
\end{lemma}

\begin{proof}
Cover a small open neighborhood of $\gamma$ by finitely many smooth flow boxes.  By passing to a smaller $\tau$ and $A$ as necessary, we may suppose that $A$ is covered by two flow boxes with union, $V$, satisfying the properties (1)-(5).
\end{proof}

\begin{notation}\label{notationattracting}
Denote the neighborhood $V$ of Lemma~\ref{holnbddefn} by $V_{\gamma}(\tau,A)$.  
\end{notation}

Notice that if $\tau=\sigma$, then $V_{\gamma}(\tau,A)$ is diffeomorphic to $A\times I$, where $I$ is a nondegenerate closed interval.  Otherwise, there is a unique smooth vertical rectangle, $R$ say, so that the result of cutting $V_{\gamma}(\tau,A)$ open along $R$, and taking the metric closure, is diffeomorphic to a solid cube.

\begin{notation} \label{notationQ}
Let $R_{\gamma}(\tau,A)$ denote any smooth vertical rectangle embedded in $V_{\gamma}(\tau,A)$ such that the result of 
cutting $V_{\gamma}(\tau,A)$ open along $R$, and taking the metric closure, is diffeomorphic to a solid cube.  When $V_{\gamma}(\tau,A)$ is not diffeomorphic to a product, $R_{\gamma}(\tau,A)$ is uniquely determined.  Let $Q_{\gamma}(\tau,A)$ denote the resulting solid cube; so

$$V_{\gamma}(\tau,A) | {R_{\gamma}(\tau,A)}=Q_{\gamma}(\tau,A).$$
\end{notation}
Note that if $\gamma$ is essential, then $Q_{\gamma}(\tau,A)$ can be viewed as a $(\tilde{\mathcal F},\tilde{\Phi})$ flow box, where $(\tilde{\mathcal F},\tilde{\Phi})$ is the lift of $(\mathcal F,\Phi)$ to the universal cover of $M$.

\begin{definition} The neighborhood $V_{\gamma}(\tau,A)$ is called the {\sl holonomy neighborhood determined by $(\tau,A)$}, and is called an {\sl attracting neighborhood} if $h_\gamma(\tau)$ is contained in the interior of $\sigma$.  
\end{definition}

\begin{definition}\label{complete attracting} Let $\mathcal F$ be a transversely oriented, $C^{\infty,0}$ foliation.  A set of holonomy neighborhoods $V_{\gamma_1}(\tau_1,A_1),...,V_{\gamma_n}(\tau_n,A_n)$ for $\mathcal F$ is {\sl spanning} if each leaf of $\mathcal F$ has nonempty intersection with the interior at least one $V_{\gamma_i}(\tau_i,A_i)$.
\end{definition}

\begin{definition}\label{Vcompatible}
Let $V$ be the union of pairwise disjoint holonomy neighborhoods $V_{\gamma_1}(\tau_1,A_1),...,V_{\gamma_n}(\tau_n,A_n)$ for $\mathcal F$.  A transversely oriented, $C^{\infty,0}$ foliation $\mathcal G$ in $M$ is called {\sl $V$-compatible with $\mathcal F$}, (or simply {\sl $V$-compatible} if $\mathcal F$ is clear from context) if each $V_{\gamma_i}(\tau_i,A_i)$ is a holonomy neighborhood for $\mathcal G$, with $V$ spanning for $\mathcal G$ if it is spanning for $\mathcal F$.  
\end{definition}

Fix a set of pairwise disjoint holonomy neighborhoods $V_{\gamma_1}(\tau_1,A_1),...,$ $ V_{\gamma_n}(\tau_n,A_n)$ for $\mathcal F$, and let $V$ denote their union.  Let $R_i=R_{\gamma_i}(\tau_i,A_i), 1\le i\le n$, and let $R$ denote the union of the $R_i$.  For each $i, 1\le i\le n$, fix a smooth open neighborhood $N_{R_i}$ of $R_i$ in $V_i$.  Choose each $N_{R_i}$ small enough so that its closure, $\overline{N_{R_i}}$, is a closed regular neighborhood of $R_i$.  Let $N_R$ denote the union of the $N_{R_i}$.

Now, given $V$, $R$ and $N_R$, we further constrain the set of foliations $\mathcal F$ (that we need to approximate by smooth contact structures) to $C^{\infty,0}$ foliations which are smooth on $N_R$.  The following lemma, applied and cited as Lemma~3.7 in \cite{KR3}, establishes that we can do this with no loss of generality.  

\begin{lemma}\label{smoothaboutR} Let $\mathcal F$ be a transversely oriented, $C^{\infty,0}$ foliation, and let $\Phi$ be a smooth flow transverse to $\mathcal F$.  Let $V$ denote the union of a set of pairwise disjoint holonomy neighborhoods for $\mathcal F$ and fix $N_R$ as above.  There is an isotopy of $M$ taking $\mathcal F$ to a $C^{\infty,0}$ foliation which is both $C^0$ close to $\mathcal F$ and smooth on $N_R$.  This isotopy may be taken to preserve $V$ and be flow compatible.
\end{lemma}

\begin{proof} This follows from the next Lemma~\ref{stronglycompatibleF}.
\end{proof}

Next we describe a preferred product parametrization on a closed set containing $V$.  In this paper, we express $S^1$ as the quotient $S^1=[-1,1]/\sim$, where $\sim$ is the equivalence relation on $[-1,1]$ which identifies $-1$ and $1$.

\begin{lemma}\label{prodnbddefn} Let $\mathcal F$ be a transversely oriented, $C^{\infty,0}$ foliation, and let $\Phi$ be a smooth flow transverse to $\mathcal F$.  Let $V$ denote the union of pairwise disjoint holonomy neighborhoods $V_i=V_{\gamma_i}(\tau_i,A_i), 1\le i\le n,$ for $\mathcal F$, and fix $N_R$ as above.  Suppose $\mathcal F$ is smooth on $N_R$.  Then for each $i, 1\le i\le n$, there is a pairwise disjoint collection of closed solid tori $P_i$ such that $P_i$ contains $V_i$ and there is a diffeomorphism $P_i\to [-1,1]\times S^1\times [-1,1]$ which satisfies the following: 
\begin{enumerate}
\item the flow segments $\Phi\cap P_i$ are identified with the segments $\{(x,y)\}\times [-1,1]$,
\item $A_i$ is identified with $[-1,1]\times S^1\times \{0\}$,
\item $\gamma_i$ is identified with $\{0\}\times S^1\times \{0\}$, 
\item $R_i$ is identified with $[-1,1]\times \{1\sim -1\}\times [-1,1]$, and
\item the restriction of the diffeomorphism to $N_{R_i}$ maps leaves of $\mathcal F$ to horizontal level surfaces $D_z = D\times \{z\}$, where $ D=[-1,1]\times (( [1/2,1] \cup [-1,-1/2])/\sim)$.
\end{enumerate}
\end{lemma}

\begin{proof}
Since $V_i$ is homeomorphic to a solid torus, it is contained in a solid torus which is diffeomorphic to $[-1,1]\times S^1\times [-1,1]$, where the diffeomorphism can be chosen to identify $A$ with $[-1,1]\times S^1\times \{0\}$ and the flow segments $\Phi\cap P$ with the segments $\{(x,y)\}\times [-1,1]$.  Moreover, since the restriction of $\mathcal F$ to $V_i \cap N_R$ is a smooth product foliation transverse to vertical fibers, and there is a unique such up to diffeomorphism, this diffeomorphism $P_i \to [-1,1]\times S^1\times [-1,1]$ can also be chosen so that the restriction of the diffeomorphism to $N_R$ maps leaves of $\mathcal F\cap N_R$ to the horizontal level surfaces $D_z$.
\end{proof}

\begin{definition} \label{productnbddefn} Fix $V$ and $N_R$ as above.  Let $P_i$ and $P_i\to [-1,1]\times S^1\times [-1,1]$ be as given in Lemma~\ref{prodnbddefn}.  Abuse notation and use the diffeomorphism to identify $P_i$ with $[-1,1]\times S^1\times [-1,1]$.  Let $\mathcal P_i$ be the product foliation of $P_i$ with leaves $([-1,1]\times S^1)\times \{t\}$, and call such a foliated solid torus, $(P_i,\mathcal P_i)$, a {\sl product neighborhood of $(V_i;N_{R_i})$.} Letting $P$ denote the union of the $P_i$ and $\mathcal P$ denote the union of the $\mathcal P_i$, call $(P,\mathcal P)$ a {\sl product neighborhood of $(V;N_R)$.} 
\end{definition}

\begin{definition}\label{stronglycompatible}
Let $\mathcal F$ be a transversely oriented, $C^{\infty,0}$ foliation and $V$ the union of pairwise disjoint, holonomy neighborhoods $V_{\gamma_i}(\tau_i,A_i), 1\le i\le k,$ for $\mathcal F$.  Let $R$ denote the union of the $R_{\gamma_i}(\tau_i,A_i), 1\le i\le k$, and let $N_R$ be an open regular neighborhood of $R$ in $V$.  Let $(P,\mathcal P)$ be a product neighborhood of $(V;N_R)$.  The foliation $\mathcal F$ is {\sl strongly $(V,P)$ compatible} if
\begin{enumerate}
\item $\mathcal F\cap N_R=\mathcal P\cap N_R$, and 
\item in the coordinates inherited from $P$, $\mathcal F\cap V$ is a product foliation $[-1,1]\times \mathcal F_0$, where $\mathcal F_0$ is a $C^{\infty,0}$ foliation of $V\cap (\{1\}\times S^1\times [-1,1])$ (i.e., $\mathcal F\cap V$ is $x$-invariant).
\end{enumerate}
\end{definition}

Given $V$, $R$ and $N_R$, we will further constrain the set of foliations $\mathcal F$ to $C^{\infty,0}$ foliations which are strongly $(V,P)$ compatible for some choice of product neighborhood $(P,\mathcal P)$.  The following lemma, applied and cited as Lemma~3.11 in \cite{KR3}, establishes that we can do this with no loss of generality; namely, after a small perturbation of $\mathcal F$, it is possible to rechoose the diffeomorphisms $P_i\to [-1,1]\times S^1\times [-1,1]$ so that $\mathcal F=\mathcal P$ on $\overline{N_R}$ and $\mathcal F$ is invariant under translation in the first coordinate.

\begin{lemma}\label{stronglycompatibleF} Let $\mathcal F$ be a transversely oriented, $C^{\infty,0}$ foliation and let $\Phi$ be a smooth flow transverse to $\mathcal F$.  Let $V$ denote the union of a set of pairwise disjoint holonomy neighborhoods for $\mathcal F$ and fix $N_R$ as above.  There is a $C^0$ small, flow compatible, isotopy of $M$ that takes $\mathcal F$ to a $C^{\infty,0}$ foliation that is $C^0$ close to $\mathcal F$ and strongly $(V,P)$ compatible for some choice of product neighborhood $(P,\mathcal P)$ of $(V;N_R)$.  This isotopy may be taken to preserve $V$.
\end{lemma}

\begin{proof} The method of Proposition~\ref{uniqueness in a product} can be applied here.  That is, on $V$ there exists a foliation $\mathcal G$ with the smoothness, restriction, and invariance properties specified in the lemma.  Moreover, such a $\mathcal G$ can be constructed to have the same holonomy representation as $\mathcal F$ with respect to each horizontal annulus of $V$ and any vertical fiber of $\Phi \cap V$.  Thus the method of Proposition~\ref{uniqueness in a product} gives the desired isotopy.
\end{proof}

\section{Denjoy blowup}

In \cite{denjoy1,denjoy2}, Denjoy gave examples of $C^1$ foliations on $T^2$ with exceptional minimal sets.  In \cite{Di}, Dippolito generalized Denjoy's method to $C^{\infty,0}$ codimension one foliations of $n$-manifolds.  This generalized construction is commonly referred to as {\sl Denjoy blowup}, and is defined precisely as follows.

\begin{definition}\label{maindefinition} Let $L$ be a countable (finite or countably infinite) union of leaves of a $C^{k,0}$ foliation $\mathcal F$ of $M$ with $k\ge 1$, and let $\Phi$ be a smooth flow transverse to $\mathcal F$.  A $C^{k,0}$ foliation, $\mathcal F'$, is a {\sl Denjoy blowup of $\mathcal F$ along $L$} if there is an open subset $U\subset M$ and a continuous collapsing map $\pi:M \to M$ satisfying the following properties:
\begin{enumerate}

\item $\mathcal F'$ is transverse to $\Phi$,

\item there is an injective map $j:L \times I \to M$ such that $j|_{L \times (0,1)}$ is a $C^k$ immersion and $j (L\times (0,1)) = U$, 

\item for each $p\in L$, $j(\{p\} \times I)$ is contained in a flow line of $\Phi$,

\item $j(L \times \{0\})$ and $j(L \times \{1\})$ are leaves of $\mathcal F'$,

\item $\pi^{-1}(p)$ is a point if $p \notin L$ and equals $j(\{p\} \times I)$ if $p \in L$,

\item $\pi$ is $\Phi$ compatible and maps leaves of $\mathcal F'$ to leaves of $\mathcal F$, 

\item $\pi$ is $C^k$ when restricted to any leaf of $\mathcal F'$, and 

\item $\pi$ is the time one map of a $\Phi$ compatible $C^k$ isotopy $\pi_t:M \to M$.  

\end{enumerate}
If $\mathcal L$ is a $C^{k,0}$ almost horizontal foliation of $L \times I$, and the pullback of the Denjoy blowup $\mathcal F'$ to $L \times I$ is $C^k$ equivalent to $\mathcal L$, then 
$\mathcal F'$ is a Denjoy blowup of $\mathcal F$ along $L$ {\sl by $\mathcal L$.}
\end{definition}

\begin{lemma}\label{canfirstclean}
Suppose $\mathcal F$ and $\mathcal G$ are $C^{1,0}$ foliations of $M$ transverse to a common smooth flow $\Phi$.  If $\mathcal F$ and $\mathcal G$ are $\Phi$ compatible isotopic, then a Denjoy blowup of $\mathcal G$ is a Denjoy blowup of $\mathcal F$.  
\end{lemma}

\begin{proof}
Denjoy blowup is defined only up to $\Phi$ compatible isotopy, and so varying a foliation by a $\Phi$ compatible isotopy does not change its Denjoy blowup.
\end{proof}

The following result extends Dippolito's generalization of Denjoy's construction (Theorem~7, \cite{Di}) in two ways.  First, it allows for foliations which are not $C^{\infty,0}$.  Second, it shows that the resulting foliation, $\mathcal F'$, can be constructed arbitrarily $C^0$ close to $\mathcal F$.

\begin{thm}[Denjoy blowup]\label{blowup} Let $\mathcal F$ be a $C^{1,0}$ foliation in a compact 3-manifold $M$.  Suppose that $\mathcal F$ is transverse to a smooth flow $\Phi$.  Let $L$ be a countable collection of leaves of $\mathcal F$, and let $\mathcal L$ be a $C^{1,0}$ almost horizontal foliation of $L \times I$.  Then there exists $C^{\infty,0}$ $\mathcal F'$ arbitrarily $C^{0}$ close to $\mathcal F$ that is a Denjoy blowup of $\mathcal F$ along $L$ by $\mathcal L$.

Moreover, if $\mathcal F$ is $C^{\infty,0}$ and $\mathcal B$ is a $(\mathcal F,\Phi)$-flow box decomposition of $M$, with flow boxes $F_1,...,F_n$, and $L$ is disjoint from $\cup_j \partial_h F_i$, then the Denjoy blowup can be chosen to be strongly $\mathcal B$ compatible in the following sense: the restriction of $\mathcal F'$ to each $F_i$ is the Denjoy blowup of the restriction of $\mathcal F$ to $F_i$.
\end{thm}

By Theorem~\ref{calegarikr}, a $C^{k,0}$ foliation on a compact 3-manifold can be isotoped by a $C^0$ small $\Phi$ compatible isotopy to a $C^0$ close, $C^{\infty,0}$ foliation.  Thus, by Lemma~\ref{canfirstclean}, it is enough to prove the theorem with the assumption that $\mathcal F$ is $C^{\infty,0}$ and $\mathcal L$ is $C^{k,0}$.  While it may be possible to conclude $C^0$ proximity of $\mathcal F'$ from Dippolito's original proof, the method of flow box decompositions gives a direct and elementary proof.

We first describe the Denjoy blowup of a strictly horizontal (and therefore smooth) foliation $\mathcal F$ of a single flow box.  

\begin{lemma} \label{one flow box} Let $\mathcal F$ be a strictly horizontal foliation of a $C^{\infty}$ flow box $F=D\times I$.  Let $L$ be a countable union of leaves of $\mathcal F$, and let $\mathcal L$ be a $C^{k,0}$ almost horizontal foliation of $L \times I$, for some $k\ge 1$.  Then there exists a $C^{\infty,0}$ foliation $\mathcal F'$ arbitrarily $C^{0}$ close to $\mathcal F$ that is a Denjoy blowup of $\mathcal F$ along $L$ by $\mathcal L$.

Moreover, given finitely many leaves $D\times \{t_j\}$ of $\mathcal F$ that are disjoint from $L$, $\mathcal F'$ can be chosen so that the restriction of $\pi$ to each $D\times \{t_j\}$ is the identity map.
\end{lemma}

\begin{proof} In this case, $\Phi$ is a flow along the vertical segments $\{x\}\times I$.  

Let $D_t=D\times \{t\}$, and let the components of $L$ be the leaves $D_{z_i}$ for some set of points $z_i\in (0,1)$, $i\in\mathcal A$.  

\begin{figure}[htbp] 
\centering
\includegraphics[width=4in]{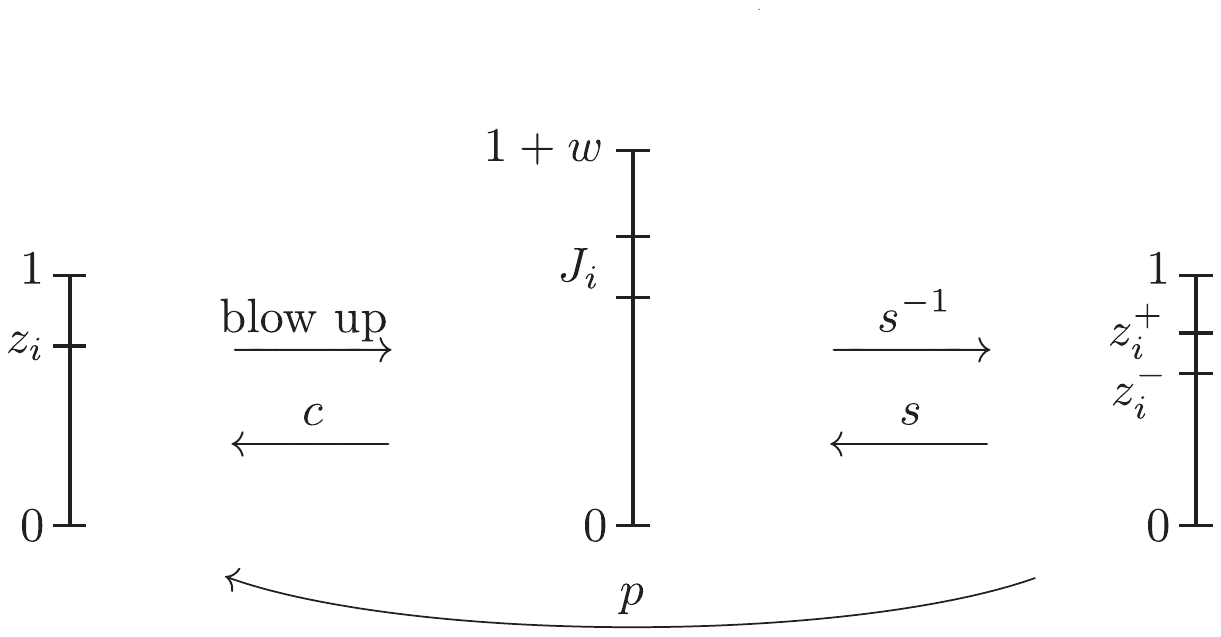} 
\caption{}
\label{blowupfig}
\end{figure}

We begin by describing the Denjoy blowup of $I$ along the points $z_i$.  Let $w_i$ denote a summable sequence of positive numbers, with sum $w=\Sigma_i w_i$.  Cut $I$ at each $z_i$ and insert an interval $J_i$ of length $w_i$.  The result is a new interval of length $1+w$.  The left inverse of this operation is a Cantor function; denote this Cantor function by $c:[0,1+w]\to [0,1]$.  Let $p:[0,1]\to [0,1]$ denote the function obtained by composing the function $c$ with the linear scaling $s:[0,1]\to [0,1+w]$; so $p=c\circ s$.  These function are illustrated in Figure~\ref{blowupfig}.

Let $[z_i^-,z_i^+]=s^{-1}(J_i)$, and notice that $z_i^+ - z_i^-=w_i/(1+w)$.  Set $C=I\setminus \sqcup_i (z_i^-,z_i^+)$, and let $\Lambda$ be the strictly horizontal lamination with leaves $D_t, t\in C$.  Let $U_i=D\times (z_i^-,z_i^+)$, and $U=\sqcup_i U_i$, the open set $F\setminus \Lambda$.

Let $\pi=id\times p:F \to F$, where $id$ is the identity map on $D$.  In particular, $\pi$ takes each flow segment $(x\times I)\cap U_i$ to the point $(x,z_i)$.

Fix $i\in \mathcal A$, and let $f_i:I\to [z_i^-,z_i^+]$ be the linear diffeomorphism.  Define $j_i=id\times f_i:(D_{z_i}\times I)\to U_i$, and define $j=\cup j_i:L\times I\to U$ to be the map that restricts to $j_i$ on $D_{z_i}\times I$.  

Let $\mathcal F'$ denote the foliation obtained by taking the union of $\Lambda$ with $j(\mathcal L)$.  Now fix $\epsilon>0$.  Since $f_i' =w_i/(1+w) < w_i$, and $\mathcal F$ is strictly horizontal, the $w_i$ can be chosen so that $\mathcal F'$ is $\epsilon$ $C^0$ close to $\mathcal F$.

Properties (1)--(7) of Definition~\ref{maindefinition} then follow immediately.  The isotopy $\pi_t$ of Property (8) is given by the straight line, $\Phi$ compatible, isotopy from the identity map to $\pi$.  By Theorem~\ref{calegarikr}, we may isotope the resulting $C^{1,0}$ Denjoy blowup to a $C^0$ close $C^{\infty,0}$ Denjoy blowup.

Finally, if $D\times \{t_j\}$ is a listing of finitely many leaves of $\mathcal F$ that are disjoint from $L$, cut $F$ open along along each $D\times \{t_j\}$, and perform Denjoy blowup, as just described, on each resulting flow box.
\end{proof} 

Hence, Theorem~\ref{blowup} holds for a strictly horizontal (and therefore smooth) foliation $\mathcal F$ of a single flow box.  

\begin{proof}[Proof of Theorem~\ref{blowup}] By Theorem~\ref{calegarikr}, $\mathcal F$ is $\Phi$ compatible isotopic to a $C^0$ close $C^{\infty,0}$ foliation.  By Lemma~\ref{canfirstclean} therefore, we may restrict attention to the case that $\mathcal F$ is $C^{\infty,0}$.  

Let $\mathcal B$ be a smooth $(\mathcal F,\Phi)$ flow box decomposition of $M$.  Let $F_1,...,F_n$ be a listing of the flow boxes of $\mathcal B$.  Choose $\mathcal B$ so that $\cup_i \partial_h F_i$ is disjoint from $L$.

Let $\sigma_i,1\le i\le n$, be a listing of the maximal vertical faces of $\mathcal B$, and let $$(N,N_v,N(\sigma_1),...,N(\sigma_n))$$ be a regular neighborhood structure for $\mathcal B$.  By Corollary~\ref{smoothtransversal2} and Lemma~\ref{canfirstclean}, it suffices to further restrict attention to the case that $\mathcal F$ is $C^{\infty,0}$, and smooth when restricted to $N_v$.

We will describe a $C^0$ close Denjoy splitting $\mathcal F'$ by considering first $N_v$, then the union $\cup_i N(\sigma_i)$, and finally the flow box interiors forming the complement of $N$.

Recall that $N_v$ is a union of flow boxes, $B_j=D_j\times I$, satisfying conditions (3) and (4) of Definition~\ref{rnd defn}.  Rechoose the $B_j$, if necessary, so that $D_j\times (0,1)$ has empty intersection with $\cup_i \partial_h F_i$.  This can be achieved by cutting each $B_j$ open along any horizontal level that has nonempty intersection with $\cup_i \partial_h F_i$.

Let $B=D\times I$ be a $B_j$ that has nonempty intersection with $L$.  Since the restriction of $\mathcal F$ to $B$ is smooth, there is a smooth parametrization $(\overline{x},z)$ of $B$ such that the restriction of $\mathcal F $ to $B$ is strictly horizontal and the restriction of $\Phi$ to $B$ has flow lines the vertical line segments $\overline{x}\times I$.  
By Lemma~\ref{one flow box}, therefore, there is a Denjoy blowup $\mathcal F'_B$ of the restriction of $\mathcal F$ to $B$, and hence functions $\pi_B:B\to B$ and $j_B:L_B\times I\to B$ satisfying the conditions of Definition~\ref{maindefinition}.  

Repeat this process for each $B_j$ that has nonempty intersection with $L$.  And let $\mathcal F'_{B_j}=\mathcal F$ on the remaining $B_j$.  Thus, we get a $C^{\infty,0}$ Denjoy blowup $\mathcal F'_{N_v}$ of the restriction of $\mathcal F$ to $N_v$ that is strongly compatible with $\mathcal B$, together with functions $\pi_{N_v}:N_v\to N_v$ and $j_{N_v}:L_{N_v}\times I\to N_v$ satisfying the conditions of Definition~\ref{maindefinition}.  The foliation $\mathcal F'_ {N_v}$ can be chosen to be $C^0$ close to the restriction of $\mathcal F$ to $N_v$.  

Next, let $\sigma$ be any maximal vertical face of $\mathcal B$, and let $\tau_{\pm}$ be the vertical edges of $\sigma$.  Let $N(\tau_-)$, respectively $N(\tau_+)$, denote the component of $N_v\cap N(\sigma)$ that contains $\tau_-$, respectively $\tau_+$.  The Denjoy blowup, $\mathcal F'_{N_v}$, is defined on $N_v$, and hence on each $N(\tau_{\pm})$.  Let $L_\ell$ be a listing of the components of $L\cap N(\sigma)$.  Set $a_\ell=L_\ell\cap N(\tau_-)$, and set $b_\ell=L_\ell\cap N(\tau_+)$.  Writing $N(\sigma)=D\times I$, where $D=I\times I$, $\sigma = \{0\}\times I\times I$, and $\alpha = \{0\}\times I$.  Orient $\alpha$ so that $ \rho_{\mathcal F}(\alpha)(a_\ell)=b_\ell$.  On each $[a_\ell^-,a_\ell^+]$, define $\rho(z)=j_{N_v}\circ \rho_{\mathcal L}(\alpha)\circ {j_{N_v}^{-1}(z)}$.  On the complement of $\cup_\ell [a_\ell^-,a_\ell^+]$, where $\pi^{-1}$ is single valued, define $\rho(z)=\pi^{-1}\circ \rho_{\mathcal F}(\alpha)\circ \pi(z)$.  

By Corollary~\ref{leafsmoothingcor}, the foliation $\mathcal F'_{N_v}$ defined on $N_v$ extends to a $C^{\infty,0}$ Denjoy blowup of $\mathcal F$ on $N_v \cup N(\sigma)$ that is $\mathcal B$ compatible, $C^0$ close to $\mathcal F$ on $N(\sigma_i)$ and satisfies $\rho_{\mathcal F'}(\alpha)=\rho$.  Repeating this construction for each maximal face $\sigma_i$ extends the definition of the $C^0$ close Denjoy blowup of $\mathcal F$ to $N$.

Finally, the Denjoy blowup $\mathcal F'_{N}$ defined on $N$ extends to a strongly $\mathcal B$ compatible Denjoy blowup $\mathcal F'$, $C^0$ close to $\mathcal F$, on each flow box $F_m$ of $\mathcal B$ by damped coning, as described in Proposition~\ref{extendleafsmoothoverthreecell}.  Moreover, since at each step we are extending over a flow box $F_m$, the resulting packet, $j_{F_m}: (L \times I) \cap F_m \to F_m$, of inserted leaves will be $C^{k,0}$ equivalent to $\mathcal L$ on $F_m$.  Similarly, the extension of the collapsing map $\pi:N \to N$ to $F_m$ is uniquely determined by the properties that it maps leaves of $\mathcal F'$ to leaves of $\mathcal F$ and maps each $I$ fiber to itself.

Since the resulting foliation $\mathcal F'$ satisfies the conditions of Definition~\ref{maindefinition} on $N$ and on each flow box of $\mathcal B$, it satisfies these conditions on $M$.  
\end{proof}

The following corollary is cited as Theorem~5.2 in \cite{KR3}.

\begin{cor} \label{Theorem5.2}
Let $\mathcal F$ be a transversely oriented, $C^{k,0}$ foliation with $k\ge 1$ that is transverse to a smooth flow $\Phi$.  Let $L$ be a countable collection of leaves of $\mathcal F$, and let $\mathcal F_1$ be a $C^{k,0}$ foliation of $L \times I$ transverse to the $I$ coordinate that contains $L \times \partial I$ as leaves.  Then there exists a $C^{k,0}$ $\mathcal F'$ arbitrarily $C^0$ close to $\mathcal F$ that is a Denjoy blowup of $\mathcal F$ along $L$, and such that the pullback of $\mathcal F'$ to $L \times I$ is equivalent to $\mathcal F_1$.

Moreover, if $V$ is the union of a set of pairwise disjoint holonomy neighborhoods for $\mathcal F$, $(W,\mathcal P)$ is a product neighborhood of $V$, and $\mathcal F$ is strongly $(V,W)$ compatible, then $\mathcal F'$ can be chosen to be strongly $(V,W)$ compatible.
\end{cor}

\begin{proof} The first paragraph of the corollary is stated as it is used in \cite{KR3}, and it follows directly from Theorem~\ref{blowup}.

For the second paragraph, it suffices to consider the case that $V$ consists of a single holonomy neighborhood.

Using the notation of \ref{notationQ}, $V$ can be cut open along $R$ into a cube $Q$.  Parametrize $Q=I^3$ so that $$N_{R}=I\times ((3/4,1]\cup [0,1/4))\times I,$$ and decompose $Q$ along $I\times \{1/2\}\times I$ into two flow boxes.  Thus $V$ is realized as a union of two flow boxes, and, by Corollary~\ref{existence of fbd}, this flow box decomposition of $V$ extends to a flow box decomposition $\mathcal B$ of $M$.  Moreover, this extension, $\mathcal B$, can be chosen so that each vertical face of $V$, except for the proper subface of $R$, is maximal.  Let $\sigma_1$ denote the maximal face $R$.  

Choose a regular neighborhood structure $(N, N_v, N(\sigma_1),...,N(\sigma_m))$ for $\mathcal B$ such that the following two properties are satisfied:
\begin{enumerate}

\item the decomposition of $N_v$ into flow boxes $B_p=D_p\times I$ satisfies: each vertical edge of $V$ appears as $0\times I\subset D_p\times I$ for some $p$, 

\item $N(\sigma_1)\cap V= N_{R}$.
\end{enumerate}

It now follows that if $\mathcal F$ is strongly $(V,W)$ compatible, then we can apply the construction of Theorem~\ref{blowup} so that the following is true:
\begin{enumerate}

\item $\mathcal F'$ is strictly horizontal in the flow boxes $B_p$ of $N_v$ that contain vertical edges of $R$, 

\item $\mathcal F'$ is strictly horizontal in $N(\sigma_1)$, and

\item in the coordinates inherited from $W$, the fixed product neighborhood of $V$, $\mathcal F'\cap V$ is $x$-invariant.
 
\end{enumerate}
Hence, if $\mathcal F$ is strongly $(V,W)$ compatible, then a $C^0$ close, $C^{k,0}$, Denjoy blowup $\mathcal F'$ can be chosen to be strongly $(V,W)$ compatible.
\end{proof}

\end{document}